\documentclass[3p]{elsarticle}
\usepackage{bm}
\usepackage[parfill]{parskip}
\usepackage[utf8]{inputenc}
\usepackage{url}
\usepackage{amsmath,amssymb,amsfonts,amsthm}
\usepackage{yhmath}
\usepackage{mathrsfs}
\usepackage{color}
\usepackage{xcolor}
\usepackage{bm}
\usepackage{graphicx}
\usepackage{subcaption}
\usepackage{lipsum}
\usepackage{amsfonts}
\usepackage{graphicx}
\usepackage{epstopdf}
\usepackage{array}
\usepackage[hidelinks]{hyperref}
\usepackage[pagewise]{lineno}
\usepackage{bbm}
\usepackage{xfrac}
\usepackage{varwidth}
\usepackage[nice]{nicefrac}

\usepackage{amsmath}
\usepackage{algorithm}
\usepackage[noend]{algpseudocode}

\ifpdf
  \DeclareGraphicsExtensions{.eps,.pdf,.png,.jpg}
\else
  \DeclareGraphicsExtensions{.eps}
\fi

\newtheorem{theorem}{Theorem}
\newtheorem{lemma}{Lemma}

\newtheorem{remark}{Remark}
\newtheorem{definition}{Definition}

\makeatletter
\newcommand{\doublehat}[1]{% 
\begingroup%
  \let\macc@kerna\z@%
  \let\macc@kernb\z@%
  \let\macc@nucleus\@empty%
  \widehat{\raisebox{.35ex}{\vphantom{\ensuremath{#1}}}\smash{\widehat{#1}}}%
\endgroup%
}
\makeatother

\makeatletter
\newsavebox\myboxA
\newsavebox\myboxB
\newlength\mylenA

\newcommand*\xoverline[2][0.75]{%
    \sbox{\myboxA}{$\m@th#2$}%
    \setbox\myboxB\null% Phantom box
    \ht\myboxB=\ht\myboxA%
    \dp\myboxB=\dp\myboxA%
    \wd\myboxB=#1\wd\myboxA% Scale phantom
    \sbox\myboxB{$\m@th\overline{\copy\myboxB}$}%  Overlined phantom
    \setlength\mylenA{\the\wd\myboxA}%   calc width diff
    \addtolength\mylenA{-\the\wd\myboxB}%
    \ifdim\wd\myboxB<\wd\myboxA%
       \rlap{\hskip 0.5\mylenA\usebox\myboxB}{\usebox\myboxA}%
    \else
        \hskip -0.5\mylenA\rlap{\usebox\myboxA}{\hskip 0.5\mylenA\usebox\myboxB}%
    \fi}
\makeatother

\usepackage{tikz}
\usetikzlibrary{chains, positioning, arrows.meta, bending, shapes.arrows}
\usetikzlibrary{positioning,shapes,arrows,calc,shapes.geometric,backgrounds,fit}
\usepackage{pifont}% http://ctan.org/pkg/pifont
\usetikzlibrary{positioning,shapes,arrows,calc,shapes.geometric,backgrounds,fit}

\tikzstyle{block} = [rectangle,draw,minimum width=2em,align=center,rounded corners, minimum height=2em,scale=1.0]
\tikzstyle{blockleft} = [rectangle,draw,minimum width=2em,align=left,rounded corners, minimum height=2em,scale=1.0]
\tikzstyle{bigblock} = [rectangle,draw,minimum width=8em,align=center,rounded corners, minimum height=4em,scale=1.0]
\tikzstyle{connect} = [draw,-latex']
\tikzstyle{decision} = [diamond, draw, 
    text width=4.5em, text badly centered, node distance=3cm, inner sep=0pt]
\tikzstyle{line} = [draw, -latex']
\tikzstyle{cloud} = [draw, ellipse,fill=red!20, node distance=3cm,
    minimum height=2em]
\tikzstyle{linenoarrow}=[draw]

\makeatletter
\let\save@mathaccent\mathaccent
\newcommand*\if@single[3]{%
  \setbox0\hbox{${\mathaccent"0362{#1}}^H$}%
  \setbox2\hbox{${\mathaccent"0362{\kern0pt#1}}^H$}%
  \ifdim\ht0=\ht2 #3\else #2\fi
  }
\newcommand*\rel@kern[1]{\kern#1\dimexpr\macc@kerna}
\newcommand*\wide@bar@[3]{%
  \begingroup
  \def\mathaccent##1##2{%
    \let\mathaccent\save@mathaccent
    \if#32 \let\macc@nucleus\first@char \fi
    \setbox\z@\hbox{$\macc@style{\macc@nucleus}_{}$}%
    \setbox\tw@\hbox{$\macc@style{\macc@nucleus}{}_{}$}%
    \dimen@\wd\tw@
    \advance\dimen@-\wd\z@
    \divide\dimen@ 3
    \@tempdima\wd\tw@
    \advance\@tempdima-\scriptspace
    \divide\@tempdima 10
    \advance\dimen@-\@tempdima
    \ifdim\dimen@>\z@ \dimen@0pt\fi
    \rel@kern{0.6}\kern-\dimen@
    \if#31
      \overline{\rel@kern{-0.6}\kern\dimen@\macc@nucleus\rel@kern{0.4}\kern\dimen@}%
      \advance\dimen@0.4\dimexpr\macc@kerna
      \let\final@kern#2%
      \ifdim\dimen@<\z@ \let\final@kern1\fi
      \if\final@kern1 \kern-\dimen@\fi
    \else
      \overline{\rel@kern{-0.6}\kern\dimen@#1}%
    \fi
  }%
  \macc@depth\@ne
  \let\math@bgroup\@empty \let\math@egroup\macc@set@skewchar
  \mathsurround\z@ \frozen@everymath{\mathgroup\macc@group\relax}%
  \macc@set@skewchar\relax
  \let\mathaccentV\macc@nested@a
  \if#31
    \macc@nested@a\relax111{#1}%
  \else
    \def\gobble@till@marker##1\endmarker{}%
    \futurelet\first@char\gobble@till@marker#1\endmarker
    \ifcat\noexpand\first@char A\else
      \def\first@char{}%
    \fi
    \macc@nested@a\relax111{\first@char}%
  \fi
  \endgroup
}
\makeatother

\newcommand{\TheTitle}{A stable multiplicative dynamical low-rank discretization for the linear Boltzmann-BGK equation} 

\date{\today}

\usepackage{amsopn}
\DeclareMathOperator{\diag}{diag}

\journal{arXiv}

\begin{document}
\begin{frontmatter}

\title{\TheTitle}

\author[adressWuerzburg]{Lena Baumann}
\author[adressInnsbruck]{Lukas Einkemmer}
\author[adressWuerzburg]{Christian Klingenberg}
\author[adressAs]{Jonas Kusch}

\address[adressWuerzburg]{Department of Mathematics, University of W\"urzburg, Emil-Fischer-Strasse 40, W\"urzburg, DE-97074, Germany, \href{mailto:lena.baumann@uni-wuerzburg.de}{lena.baumann@uni-wuerzburg.de} (Lena Baumann), \href{mailto:christian.klingenberg@uni-wuerzburg.de}{christian.klingenberg@uni-wuerzburg.de} (Christian Klingenberg) }
\address[adressInnsbruck]{Numerical Analysis and Scientific Computing, University of Innsbruck, Technikerstrasse 13, Innsbruck, A-6020, Austria, \href{mailto:lukas.einkemmer@uibk.ac.at}{lukas.einkemmer@uibk.ac.at}}
\address[adressAs]{Scientific Computing, Norwegian University of Life Sciences, Drøbakveien 31, \r{A}s, NO-1433, Norway, \href{mailto:jonas.kusch@nmbu.no}{jonas.kusch@nmbu.no}}

\begin{abstract}
The numerical method of dynamical low-rank approximation (DLRA) has recently been applied to various kinetic equations showing a significant reduction of the computational effort. In this paper, we apply this concept to the linear Boltzmann-Bhatnagar-Gross-Krook (Boltzmann-BGK) equation which due its high dimensionality is challenging to solve. Inspired by the special structure of the non-linear Boltzmann-BGK problem, we consider a multiplicative splitting of the distribution function. We propose a rank-adaptive DLRA scheme making use of the basis update \& Galerkin integrator and combine it with an additional basis augmentation to ensure numerical stability, for which an analytical proof is given and a classical hyperbolic Courant–Friedrichs–Lewy (CFL) condition is derived. This allows for a further acceleration of computational times and a better understanding of the underlying problem in finding a suitable discretization of the system. Numerical results of a series of different test examples confirm the accuracy and efficiency of the proposed method compared to the numerical solution of the full system. 
\end{abstract}

\end{frontmatter}

\section{Introduction}\label{sec1:Introduction}

Numerically solving kinetic equations usually requires immense computational and memory efforts due to the high-dimensional phase space containing all possible states of the system. The state of a kinetic system is described by a distribution function $f$ which can be interpreted as the corresponding particle density in phase space. Instead of solving one high-dimensional equation, the concept of dynamical low-rank approximation (DLRA) \cite{koch2007dynamical} allows us to split the problem into three lower dimensional subequations leading to an appropriate approximation of the solution. In particular, in a one-dimensional setting we approximate the distribution function $f(t,x,v)$, with $t \in \mathbb{R}^+$ denoting the time, $x \in D \subset \mathbb{R}$ the spatial, and $v\in \mathbb{R}$ the velocity variable, by
\begin{align*}
f(t,x,v) \approx \sum_{i,j=1}^r X_i(t,x) S_{ij}(t) V_j(t,v), 
\end{align*}
and evolve the corresponding low-rank factors in three substeps further in time. The sets $\{X_i : i=1,..,r \}$ and $\{V_j : j=1,..,r \}$ contain the orthonormal basis functions in space and in velocity, respectively, and $r$ is called the rank of this approximation. DLRA has recently gained increasing interest and has been studied in various fields including radiation transport \cite{baumann2024, Patwardhan2024, einkemmerhukusch2024,peng2020-2D, HauckSchnake2023}, radiation therapy \cite{stammerkusch2023}, plasma physics \cite{einkemmerlubich2018, einkemmerostermann2020, einkemmerjoseph2021conservative, einkemmerscalone2023}, chemical kinetics \cite{prugger2023, prugger2024} and Boltzmann type transport problems \cite{einkemmer2019weakly, einkemmerhuying2021, dingeinkemmerli2021, huwang2022}. The core idea of this method is to project the solution to a manifold of low-rank functions of the above form and constrain the solution dynamics there. Different time integrators which are able to ensure this behaviour and are robust to the presence of small singular values are available. Frequently used integrators for kinetic problems are the \textit{projector-splitting} \cite{lubich2014projector} as well as the \textit{(rank-adaptive) basis update \& Galerkin} (BUG) \cite{ceruti2022unconventional, ceruti2022rank} and the \textit{parallel} integrator \cite{ceruti2024parallel}. For the rank-adaptive BUG and the parallel integrator extensions to schemes with proven second-order robust error bounds have been derived in \cite{ceruti2024secondorderBUG,kusch2024secondorderparallel}.\\
For a large number of collisions, the solution $f$ of the Boltzmann-BGK equation stays close to the Maxwellian equilibrium distribution $M$ which in general is not a low-rank function. Inspired by \cite{einkemmerhuying2021, kormannsonnendruecker2016}, we use the multiplicative splitting $f=Mg$, for which in \cite{einkemmerhuying2021} it has been shown that $g$ is a low-rank function even if this if not the case for $f$. Hence, we derive an evolution equation for $g$ and apply the low-rank approach to this part of the distribution function. Difficulties may arise in the discretization as it is per se not clear how to treat spatial derivatives.

In this paper we propose a stable dynamical low-rank discretization for the linear Boltzmann-BGK equation. The main features of this paper are: 
\begin{itemize}
    \item \textit{A multiplicative splitting of the distribution function:} As the Maxwellian equilibrium distribution $M$ is generally not a low-rank function, we consider the multiplicative splitting $f=Mg$ and apply the low-rank ansatz to the remaining function $g$. It can be considered as a deviation from the equilibrium and is shown to be of low rank \cite{einkemmerhuying2021}.
    \item \textit{A stable numerical scheme for linear Boltzmann-BGK with rigorous mathematical proofs:} We show that a stable discretization has to be derived  carefully and compare it with an intuitive discretization that fails to guarantee numerical stability. We give a rigorous analytical proof of stability and derive a classic hyperbolic CFL condition. This enables us to choose an optimal time step size of $\Delta t = \text{CFL} \cdot \Delta x$ with $\text{CFL}$ denoting the CFL number, leading to a reduction of the computational effort.
    \item \textit{A rank-adaptive integrator:} For the low-rank scheme we use the rank-adaptive BUG integrator from \cite{ceruti2022rank}, leading to a basis augmentation in both the $K$- and $L$-step of the low-rank algorithm. Compared to the projector-splitting integrator used in \cite{einkemmer2019weakly, einkemmerhuying2021,dingeinkemmerli2021}, this allows us to determine the rank adaptively in each step avoiding the a priori choice of a certain fixed rank.
    \item \textit{A series of numerical experiments validating the derived properties:} We give a number of numerical examples that validate the derived stability while showing a significant reduction of computational and memory requirements of the low-rank scheme compared to the full order method. 
\end{itemize}

The paper is structured as follows: After the introduction in Section \ref{sec1:Introduction}, we provide background information on the linear Boltzmann-BGK equation, explain the considered multiplicative structure, and derive two possible systems of equations in Section \ref{sec2:Background}. Both systems are equivalent in the continuous setting. In Section \ref{sec3:Discretization}, we discretize in velocity and in space, before subsequently time is discretized giving two different fully discretized schemes. It is then shown in Section \ref{sec4:Stability} that a naive discretization can lead to a numerical scheme that is not von Neumann stable whereas a more careful treatment guarantees numerical stability. Section \ref{sec5:DLRAforMg} gives a brief introduction to the concept of DLRA and applies this method such that a numerically stable low-rank scheme is obtained. Numerical experiments in both $1$D and $2$D in Section \ref{sec6:NumericalResults} confirm the derived results. Section \ref{sec7:Outlook} gives a brief conclusion and outlook.  

\section{Linear Boltzmann-BGK}\label{sec2:Background}

The Boltzmann equation is a fundamental model in kinetic theory describing a gas that is not in thermodynamic equilibrium \cite{cercignani1988, struchtrup2005}. In its full formulation it makes use of the so called ``Stosszahlansatz'' leading to a collision operator for which the solution of the Boltzmann equation is demanding. To overcome this, the \textit{BGK model} \cite{bgk1954}, named after Bhatnagar, Gross and Krook, can be considered. It simplifies the collision term while maintaining the key properties of the equation. In a one-dimensional setting it reads
\begin{subequations}
\begin{align}\label{eq:linearBGKf}
\partial_t f(t,x,v) + v \partial_x f(t,x,v) &= \sigma \left(M[f](t,x,v)-f(t,x,v) \right),
\end{align}
where $f(t,x,v)$ denotes the distribution function depending on the time $t \in \mathbb{R}^+$, the spatial variable $x \in D \subset \mathbb{R}$ and the velocity variable $v \in \mathbb{R}$. The constant $\sigma$ describes the collisionality of the particles and $M[f]$ stands for the Maxwellian equilibrium distribution. It depends on the density $\rho(t,x) = \int_\mathbb{R} f(t,x,v) \mathrm{d}v$ for which we obtain an evolution equation by integrating \eqref{eq:linearBGKf} with respect to $v$. This gives
\begin{align}\label{eq:rho}
\partial_t \rho(t,x) = - \partial_x \int v f(t,x,v) \mathrm{d}v.
\end{align}
\end{subequations}

In \cite{einkemmerhuying2021} it has been shown that using the multiplicative decomposition
\begin{align}\label{eq:splittingMg}
f(t,x,v)=M[f](t,x,v)g(t,x,v)
\end{align}
is advantageous as $g$ is low-rank even if this is not the case for the Maxwellian (which is not true for the classic additive micro-macro decomposition). In order to reduce computational and memory costs a dynamical low-rank approach has then been applied in \cite{einkemmerhuying2021} to treat the resulting evolution equations for $g$.

In this work, we consider an isothermal Maxwellian without drift, i.e.
\begin{align*}
M[f](t,x,v) = \frac{\rho(t,x)}{\sqrt{2\pi}} \exp{(-v^2/2)}
\end{align*}
with strictly positive density $\rho(t,x) >0$. This results in a linear model, which we call the \textit{linear Boltzmann-BGK equation}. It has been extensively studied in the PDE community (see, e.g., \cite{evans2021, evans2020, achleitner2016}) as well as from a numerical point of view \cite{ bessemoulin2020}. The linearity of the model allows for a rigorous theoretical stability analysis in the sense of von Neumann, which, to our knowledge, is only available in the linear case. In this paper, we provide such stability considerations in the context of dynamical low-rank simulation using a multiplicative decomposition. The stability analysis for the simplified problem provides insight into the numerical schemes that have been used in the literature \cite{einkemmerhuying2021} for the non-linear Boltzmann-BGK equation. In particular, it explains why such multiplicative schemes need to take relatively small time step sizes even though the collision operator is treated implicitly. In addition, we expect that ideas for the design of stable numerical schemes for the Boltzmann-BGK equation can be transferred from this provably stable linear problem to more complicated non-linear setups.

We insert the multiplicative approach \eqref{eq:splittingMg} into \eqref{eq:linearBGKf} and \eqref{eq:rho} and obtain
\begin{subequations}\label{eqs-naive:Mg-both}
\begin{align}
\partial_t g(t,x,v) &= - v \partial_x  g(t,x,v) + \sigma \left(1-g(t,x,v)\right) - \frac{g(t,x,v)}{\rho(t,x)} \partial_t \rho(t,x) - v \frac{g(t,x,v)}{\rho(t,x)} \partial_x \rho(t,x),\label{eqs-naive:Mg-g}\\
\partial_t \rho(t,x) &= - \frac{1}{\sqrt{2\pi}} \partial_x \int \rho(t,x) g(t,x,v) v e^{-v^2/2} \mathrm{d}v.\label{eqs-naive:Mg-rho}
\end{align}
\end{subequations}
This set of equations is called the \textit{advection form} of the multiplicative system. It corresponds to the way the equations are treated in \cite{einkemmerhuying2021}. We can rewrite equation \eqref{eqs-naive:Mg-g} into a \textit{conservative form}, leading to the system 
\begin{subequations}\label{eqs:Mg-both}
\begin{align}
\partial_t g(t,x,v) &= - \frac{v}{\rho(t,x)} \partial_x \left(\rho(t,x) g(t,x,v)\right) + \sigma \left(1-g(t,x,v)\right) - \frac{g(t,x,v)}{\rho(t,x)} \partial_t \rho(t,x),\label{eqs:Mg-g}\\
\partial_t \rho(t,x) &= - \frac{1}{\sqrt{2\pi}} \partial_x \int \rho(t,x) g(t,x,v) v e^{-v^2/2} \mathrm{d}v.\label{eqs:Mg-rho}
\end{align}
\end{subequations}
Note that for both systems we omit initial and boundary conditions for now. It is a challenging task to construct a suitable numerical scheme as it is per se not clear how to treat the spatial derivative in the transport part of the first equations and which of both systems to prefer. Further, the potentially stiff collision term requires an implicit time discretization. In addition, the consideration of higher dimensionality occurring in practical applications leads to prohibitive numerical costs. To overcome this last problem, we make use of the numerical reduced order method of dynamical low-rank approximation after having derived a stable discretization.

\section{Discretization of the $Mg$ system}\label{sec3:Discretization}

In this section we give a full discretization of both versions \eqref{eqs-naive:Mg-both} and \eqref{eqs:Mg-both} of the $Mg$ system. We start with discretizing equations \eqref{eqs-naive:Mg-both} and \eqref{eqs:Mg-both} in velocity and space before a time discretization is presented in the next subsection.

\subsection{Discretization in velocity and in space}

For the discretization in the velocity space we use a nodal approach and prescribe a certain number of grid points $N_v \in \mathbb{N}$. Due to the special structure of \eqref{eqs-naive:Mg-rho} and \eqref{eqs:Mg-rho} we use a Gauss-Hermite quadrature providing the quadrature nodes $v_1,...,v_{N_v}$ and weights $\omega_1,...,\omega_{N_v}$ enabling us to approximate integrals as
\begin{align*}
    \int_\mathbb{R} e^{-v^2} g(t,x,v) \mathrm{d}v \approx \sum_{k=1}^{N_v} \omega_k g(t,x,v_k).
\end{align*}
For the discretization of the spatial domain $D \subset \mathbb{R}$ we take $N_x \in \mathbb{N}$ grid points and choose a grid $x_1,...,x_{N_x}$ with equidistant spacing $\Delta x = \frac{1}{N_x}$. We approximate $x$-dependent quantities by
\begin{align*}
    \rho_j(t) \approx \rho(t,x_j) \qquad \text{ and } \qquad g_{jk}(t) \approx g(t,x_j,v_k).
\end{align*} 
Spatial derivatives $\partial_x$ are approximated by the tridiagonal stencil matrices $\mathbf{D}^x \in \mathbb{R}^{N_x \times N_x}$ corresponding to a first-order central differencing scheme. Further, a tridiagonal second-order  central differencing stabilization matrix $\mathbf{D}^{xx} \in \mathbb{R}^{N_x \times N_x} $ approximating $ \partial_{xx}$ is added. Their entries are defined as
\begin{align*}
D_{j,j\pm 1}^{x}= \frac{\pm 1}{2\Delta x}\;,\qquad D_{j,j}^{xx}= - \frac{2}{\left(\Delta x\right)^2}\;, \quad D_{j,j\pm 1}^{xx}= \frac{1}{\left(\Delta x\right)^2}\;,
\end{align*}
whereas all other entries are set to zero. Note that from now on we assume periodic boundary conditions. For this reason we set 
\begin{align*}
D_{1,N_x}^{x} &= \frac{-1}{2\Delta x}
\;,\quad D_{N_x,1}^{x} = \frac{1}{2\Delta x}\;,\\
D_{1,N_x}^{xx} &= D_{N_x,1}^{xx} = \frac{1}{\left( \Delta x \right)^2}.
\end{align*}
Similar to the proof in \cite{baumann2024}, one can show that the stencil matrices $\mathbf{D}^x$ and $\mathbf{D}^{xx}$ fulfill the following properties:

\begin{lemma}[Summation by parts]\label{lemma:stencilmatrices}
Let $y,z \in \mathbb{R}^{N_x}$ with indices $i,j=1,...,N_x$. %In addition, we set $y_0 = y_{n_x}$ and $y_{n+1} = y_1$, for $z$ respectively, due to the periodic boundary conditions.
Then it holds
\begin{align*}
\sum_{i,j=1}^{N_x} y_j D_{ji}^x z_i = - \sum_{i,j=1}^{N_x} z_j D_{ji}^{x} y_i\;, \qquad \sum_{i,j=1}^{N_x} z_j D_{ji}^{x} z_i = 0 \;, \qquad \sum_{i,j=1}^{N_x} y_j D_{ji}^{xx} z_i = \sum_{i,j=1}^{N_x} z_j D_{ji}^{xx} y_i.
\end{align*}
Moreover, let $\mathbf D^{+}\in\mathbb{R}^{N_x \times N_x}$ be defined as
\begin{align*}
D_{j,j}^{+}= \frac{- 1}{\Delta x}\;,\qquad D_{j,j + 1}^{+}= \frac{ 1}{\Delta x}\;.
\end{align*}
Then, $\sum_{i,j =1}^{N_x} z_j D_{ji}^{xx} z_i = -  \sum_{j=1}^{N_x} \left(\sum_{i=1}^{N_x} D_{ji}^+ z_i\right)^2$.
\end{lemma}

We insert the proposed velocity and space discretization into the advection form \eqref{eqs-naive:Mg-both} and add a stabilizing second-order term for $\partial_x g$. This corresponds to the method used in \cite{einkemmerhuying2021} for the non-linear isothermal Boltzmann-BGK equation and leads to the semi-discrete time-continuous system
\begin{subequations}\label{eqs-naive:Mg-both-semidiscrete}
\begin{align}
\dot g_{jk}(t) =& - \sum_{i=1}^{N_x} 
D_{ji}^x g_{ik}(t) v_k + \frac{\Delta x}{2} \sum_{i=1}^{N_x} D_{ji}^{xx} g_{ik}(t) |v_k|+ \sigma \left(1-g_{jk}(t)\right) - \frac{g_{jk}(t)}{\rho_j(t)} \dot \rho_j(t)  -\sum_{i=1}^{N_x} \frac{g_{jk}(t)}{\rho_j(t)} D_{ji}^x \rho_i(t)  v_k,\label{eqs-naive:Mg-g-semidiscrete}\\
\dot \rho_j(t) =& - \frac{1}{\sqrt{2\pi}} \sum_{i=1}^{N_x} \sum_{k=1}^{N_v} D_{ji}^x \rho_i(t) g_{ik}(t) v_k \omega_k e^{v_k^2/2} + \frac{\Delta x}{2\sqrt{2\pi}} \sum_{i=1}^{N_x} \sum_{k=1}^{N_v} D_{ji}^{xx} \rho_i(t) g_{ik}(t) |v_k| \omega_k e^{v_k^2/2}.\label{eqs-naive:Mg-rho-semidiscrete}
\end{align}
\end{subequations}
For the conservative form \eqref{eqs:Mg-both} the second-order stabilization term is applied to $\partial_x (\rho g)$. We obtain the semi-discrete system
\begin{subequations}\label{eqs:Mg-both-semidiscrete}
\begin{align}
\dot g_{jk}(t) =& - \sum_{i=1}^{N_x} \frac{1}{\rho_j(t)} D_{ji}^x \rho_i(t) g_{ik}(t) v_k + \frac{\Delta x}{2} \sum_{i=1}^{N_x} \frac{1}{\rho_j(t)} D_{ji}^{xx} \rho_i(t) g_{ik}(t) |v_k| + \sigma \left(1-g_{jk}(t)\right) - \frac{g_{jk}(t)}{\rho_j(t)} \dot \rho_j(t),\label{eqs:Mg-g-semidiscrete}\\
\dot \rho_j(t) =& - \frac{1}{\sqrt{2\pi}} \sum_{i=1}^{N_x} \sum_{k=1}^{N_v} D_{ji}^x \rho_i(t) g_{ik}(t) v_k \omega_k e^{v_k^2/2} + \frac{\Delta x}{2 \sqrt{2\pi}} \sum_{i=1}^{N_x} \sum_{k=1}^{N_v} D_{ji}^{xx} \rho_i(t) g_{ik}(t) |v_k| \omega_k e^{v_k^2/2}.\label{eqs:Mg-rho-semidiscrete}
\end{align}
\end{subequations}

\subsection{Time discretization}

The time discretization of both systems has to be derived carefully to ensure numerical stability. We start with the advection form \eqref{eqs-naive:Mg-both-semidiscrete} and perform an explicit Euler step for the transport part in \eqref{eqs-naive:Mg-g-semidiscrete} as well as in \eqref{eqs-naive:Mg-rho-semidiscrete}. The potentially stiff collision term is treated implicitly. For approximating the time derivative $\partial_t \rho$ the corresponding difference quotient is used. We obtain the fully discrete scheme
\begin{subequations}
\begin{align}
g_{jk}^{n+1} =& \ g_{jk}^n -\Delta t \sum_{i=1}^{N_x} D_{ji}^x g_{ik}^n  v_k + \Delta t \frac{\Delta x}{2} \sum_{i=1}^{N_x}  D_{ji}^{xx} g_{ik}^n |v_k|\label{eqs:gEH}\\
&+ \sigma \Delta t \left(1- g_{jk}^{n+1}\right) - 
\frac{g_{jk}^{n+1}}{\rho_j^n} \left(\rho_j^{n+1}-\rho_j^n \right)- \Delta t \frac{g_{jk}^n }{\rho_j^n} \sum_{i=1}^{N_x} D_{ji}^x \rho_i^n  v_k,\nonumber\\
\rho_j^{n+1} =& \ \rho_j^n - \Delta t \frac{1}{\sqrt{2\pi}} \sum_{i=1}^{N_x} \sum_{k=1}^{N_v} D_{ji}^x \rho_i^n g_{ik}^n v_k \omega_k e^{v_k^2/2} + \Delta t \frac{\Delta x}{2 \sqrt{2\pi}} \sum_{i=1}^{N_x} \sum_{k=1}^{N_v} D_{ji}^{xx} \rho_i^n g_{ik}^n |v_k| \omega_k e^{v_k^2/2}.\label{eqs:rhoEH}
\end{align}\label{eqs:EH-both}   
\end{subequations}

For the conservative form \eqref{eqs:Mg-both-semidiscrete} we again perform an explicit Euler step for the transport part in \eqref{eqs:Mg-g-semidiscrete} as well as in \eqref{eqs:Mg-rho-semidiscrete}. The collision term is treated implicitly and a factor $\frac{\rho^{n+1}}{\rho^n}$ coming from the analysis is added. As before, the time derivative $\partial_t \rho$ is approximated by its difference quotient. This leads to the fully discretized equations
\begin{subequations}
\begin{align}
g_{jk}^{n+1} =& \ g_{jk}^n -\Delta t \sum_{i=1}^{N_x} \frac{1}{\rho_j^n} D_{ji}^x \rho_i^n g_{ik}^n v_k + \Delta t \frac{\Delta x}{2} \sum_{i=1}^{N_x} \frac{1}{\rho_j^n} D_{ji}^{xx} \rho_i^n g_{ik}^n |v_k|\label{eqs:gKB}\\
&+ \sigma \Delta t \frac{\rho_j^{n+1}}{\rho_j^n} \left(1 -g_{jk}^{n+1}\right) - \frac{g_{jk}^{n+1}}{\rho_j^n} \left(\rho_j^{n+1}-\rho_j^n \right) ,\nonumber\\
\rho_j^{n+1} =& \ \rho_j^n - \Delta t \frac{1}{\sqrt{2\pi}} \sum_{i=1}^{N_x} \sum_{k=1}^{N_v} D_{ji}^x \rho_i^n g_{ik}^n v_k \omega_k e^{v_k^2/2} + \Delta t \frac{\Delta x}{2 \sqrt{2\pi}} \sum_{i=1}^{N_x} \sum_{k=1}^{N_v} D_{ji}^{xx} \rho_i^n g_{ik}^n |v_k| \omega_k e^{v_k^2/2}.\label{eqs:rhoKB}
\end{align}\label{eqs:KB-both}   
\end{subequations}

Note that the discretizations for $\rho$ given in \eqref{eqs:rhoEH} and \eqref{eqs:rhoKB} are exactly the same. The main difference between the naive discretization of the advection form \eqref{eqs:EH-both} and the proposed scheme \eqref{eqs:KB-both} is the stabilization of $\partial_x (\rho g)$ in \eqref{eqs:gKB}, opposed to a stabilization of $\partial_x g$ as done in \eqref{eqs:gEH}. This leads to an upwind type discretization for \eqref{eqs:gKB}, whereas \eqref{eqs:gEH} corresponds to a centered finite difference scheme. Moreover, the schemes are distinguished by an additional factor $\frac{\rho^{n+1}}{\rho^n}$ in the collision term of \eqref{eqs:gKB} which is essential for the proof of numerical stability.

In the fully discrete setting, we make use of the following notations.
\begin{definition}[Fully discrete solution and Maxwellian]
The full solution $f$ of the linear Boltzmann-BGK equation in the fully discrete setting at time $t_n$ is given by $\mathbf{f}^n = (f_{jk}^n) \in \mathbb{R}^{N_x \times N_v}$ with entries 
\begin{align*}
f_{jk}^n = \frac{1}{\sqrt{2\pi}} \rho_j^n g_{jk}^n e^{-v_k^2/2}. 
\end{align*}
For the fully discrete Maxwellian $\mathbf{M}^n = (M_{jk}^n) \in \mathbb{R}^{N_x \times N_v}$ at time $t_n$, we have $M_{jk}^n = \frac{1}{\sqrt{2\pi}} \rho_j^n e^{-v_k^2/2}$. 
\end{definition}

\section{Numerical stability}\label{sec4:Stability}

Although the derivation of the equations in \eqref{eqs:EH-both} and \eqref{eqs:KB-both} is similar, both systems differ drastically in terms of numerical stability. In this section, both fully discretized schemes presented are compared.

\subsection{Naive discretization}

We begin with the naive discretization given in \eqref{eqs:EH-both} that is comparable to the one chosen in \cite{einkemmerhuying2021} in the sense that the advection form of the multiplicative splitting is used. In \cite{einkemmerhuying2021}, numerical experiments are given but no explicit stability analysis is conducted. In the following, we give an example which shows that numerical stability in the sense of von Neumann can not be guaranteed. 

\begin{theorem}\label{theorem:EH-Counterexample}
There exist initial values $\mathbf{g}^n = \left( g_{jk}^n \right) \in \mathbb{R}^{N_x \times N_v}$ and $\boldsymbol{\rho}^n = \left( \rho_j^n \right) \in \mathbb{R}^{N_x}$ such that the numerical scheme proposed in \eqref{eqs:EH-both} for $\sigma =0$ is not von Neumann stable. %in the sense that 
%\begin{align*}
%\Vert \mathbf{M}^{n+1} \cdot \mathbf{g}^{n+1} \Vert_{\mathscr{H}}^2 > \Vert \mathbf{M}^n \cdot \mathbf{g}^n \Vert_{\mathscr{H}}^2,
%\end{align*}
%where $\mathbf{M}^{n,n+1} = \left(M_{jk}^{n,n+1} \right)$ with $M_{jk}^{n,n+1} = \frac{1}{\sqrt{2\pi}} \rho_j^{n,n+1} e^{-v_k^2/2}$ and $\mathbf{g}^{n,n+1} = \left(g_{jk}^{n,n+1} \right)$.
\end{theorem}
\begin{proof}
Let us assume a solution $g_{jk}^n$ that is constant in space and velocity, e.g. $g_{jk}^n \equiv 1 $. For this solution the terms containing $\mathbf{D}^x \mathbf{g}^n$ and $\mathbf{D}^{xx} \mathbf{g}^n$ are zero. Let us further assume that there is no collisionality, i.e. $\sigma=0$. We insert this information into \eqref{eqs:gEH} and get
\begin{align*}
g_{jk}^{n+1} =& \ 1 - 
\frac{g_{jk}^{n+1}}{\rho_j^n} \left(\rho_j^{n+1}-\rho_j^n \right) - \Delta t \frac{1}{\rho_j^n} \sum_{i=1}^{N_x} \left(D_{ji}^x \rho_i^n \right) v_k.
\end{align*}
 
After rearranging, we have that 
\begin{align*}
\rho_j^{n+1} g_{jk}^{n+1} = \rho_j^n - \Delta t  \sum_{i=1}^{N_x} \left(D_{ji}^x \rho_i^n \right) v_k.
\end{align*}
Multiplication with $\frac{1}{\sqrt{2\pi}} e^{-v_k^2/2}$ then leads to
\begin{align}\label{eq:EH-Counterexample}
f_{jk}^{n+1} = f_{jk}^n - \Delta t  \sum_{i=1}^{N_x} D_{ji}^x f_{ik}^n v_k.
\end{align}
This corresponds to a discretization of $\partial_t f + v \partial_x f = 0$ with explicit Euler in time and centered finite differences in space for which it is well-known that it is not von Neumann stable \cite{hairerwanner1996, leveque2002}.
\end{proof}
Indeed, one can show that the discretization given in \eqref{eq:EH-Counterexample} is not von Neumann stable, but stable for relatively small time step sizes \cite{leveque2002}. This matches our numerical insights from \cite{einkemmerhuying2021}, where the space discretization is comparable to \eqref{eqs:EH-both} and small time step sizes are required.

\subsection{Stable discretization}

Having seen that for a certain choice of the initial values the system of equations \eqref{eqs:EH-both} is not von Neumann stable, we now consider equations \eqref{eqs:KB-both} in terms of numerical stability. We observe that the advection terms are treated explicitly, whereas the collision term is treated implicitly, leading to a removal of the potential stiffness caused by a large number of collisions. We seek a rigorous proof of stability under a classic hyperbolic CFL condition that will be derived in the following norm.
\begin{definition}[Stability norm]
For $\mathbf{f}^n = (f_{jk}^n) \in \mathbb{R}^{N_x \times N_v}$, the $\mathscr{H}$-norm shall be defined as 
\begin{align*}
    \Vert \mathbf{f} ^n \Vert_{\mathscr{H}}^2 = \sqrt{2\pi} \sum_{j=1}^{N_x} \sum_{k=1}^{N_v} \left(f_{jk}^n \right)^2 \omega_k e^{3v_k^2/2}. 
\end{align*}
This corresponds to a Frobenius norm $\Vert \cdot \Vert_F$ with weights $\sqrt{2\pi} \omega_k e^{3v_k^2/2}$.  
\end{definition}
The choice of this norm is inspired by the analysis in \cite{achleitner2016}, where hypocoercivity for the linear Boltzmann-BGK equation is shown. Different from the considerations in \cite{achleitner2016}, we use a fully discrete analogue to the considered weighted $L^2$-norm that also takes the Gauss-Hermite quadrature  into account. Note that the factor $\sqrt{2\pi}$ does not affect the stability but is added for consistency.

At each time step the fully discrete distribution function $f$ and the fully discrete density $\rho$ are required to fulfill the relation
\begin{align}\label{4:identityrhof}
\rho_j^n = \sum_{k=1}^{N_v} f_{jk}^n \omega_k e^{v_k^2} \qquad \text{ for all } n \in \mathbb{N},
\end{align}
which is the discrete counterpart of the identity $\rho = \int_{\mathbb{R}} f \mathrm{d}v$. Inserting the multiplicative splitting approach, the latter continuous identity can be rewritten as $1 = \frac{1}{\sqrt{2\pi}} \int_{\mathbb{R}} g e^{-v^2/2} \mathrm{d}v$. Hence, rewriting relation \eqref{4:identityrhof} in the same way, yields the equivalent formulation
\begin{align}\label{4:identityrhog}
1 = \frac{1}{\sqrt{2\pi}} \sum_{k=1}^{N_v} g_{jk}^n \omega_k e^{v_k^2/2} \qquad \text{ for all } n \in \mathbb{N}.
\end{align}
We are able to show that the equality given in \eqref{4:identityrhog} holds for the conservative equations
\eqref{eqs:KB-both} under a suitable choice of the initial condition.

\begin{lemma}\label{lemma:fullydiscrete-rho}
Let us assume that the initial condition for $g$ satisfies
\begin{align*}
1 = \frac{1}{\sqrt{2\pi}} \sum_{k=1}^{N_v} g_{jk}^0 \omega_k e^{v_k^2/2} \qquad \text{ for all } j \in \{ 1,...,N_x\}.
\end{align*}
Then, for all $n \in \mathbb{N}$, the equality given in \eqref{4:identityrhog} holds.
\end{lemma}
\begin{proof}
The proof follows by induction. For the induction assumption let us assume that the relation $1 = \frac{1}{\sqrt{2\pi}} \sum_{k=1}^{N_v} g_{jk}^n \omega_k e^{v_k^2/2}$ holds for one $n \in \mathbb{N}$. For the induction step we begin with equation \eqref{eqs:gKB}, put the terms containing $g_{jk}^{n+1}$ to the left-hand side and multiply with $\rho_j^{n+1}$. This results in
\begin{align*}
\left(1+ \sigma \Delta t\right) \rho_j^{n+1} g_{jk}^{n+1} =& \ \rho_j^n g_{jk}^n -\Delta t \sum_{i=1}^{N_x} D_{ji}^x \rho_i^n g_{ik}^n v_k + \Delta t \frac{\Delta x}{2} \sum_{i=1}^{N_x} D_{ji}^{xx} \rho_i^n g_{ik}^n  |v_k| + \sigma \Delta t \rho_j^{n+1}.
\end{align*}
Multiplication with $\frac{1}{\sqrt{2\pi}} \omega_k e^{v_k^2/2}$ and summation over $k$ leads to
\begin{align*}
\left(1+ \sigma \Delta t\right) \rho_j^{n+1} \frac{1}{\sqrt{2\pi}} \sum_{k=1}^{N_v} g_{jk}^{n+1} \omega_k e^{v_k^2/2} =& \ \rho_j^n \frac{1}{\sqrt{2\pi}} \sum_{k=1}^{N_v} g_{jk}^n \omega_k e^{v_k^2/2} -\Delta t \frac{1}{\sqrt{2\pi}} \sum_{i=1}^{N_x} \sum_{k=1}^{N_v} D_{ji}^x \rho_i^n g_{ik}^n v_k \omega_k e^{v_k^2/2}\\
&+ \Delta t \frac{\Delta x}{2 \sqrt{2\pi}} \sum_{i=1}^{N_x} \sum_{k=1}^{N_v} D_{ji}^{xx} \rho_i^n g_{ik}^n |v_k| \omega_k e^{v_k^2/2}\\
&+ \sigma \Delta t \rho_j^{n+1} \frac{1}{\sqrt{2\pi}} \sum_{k=1}^{N_v} \omega_k e^{v_k^2/2}. 
\end{align*}
We insert the induction assumption as well as $\frac{1}{\sqrt{2\pi}} \sum_{k=1}^{N_v} \omega_k e^{v_k^2/2} = 1$. Then, together with equation \eqref{eqs:rhoKB}, this establishes
\begin{align*}
\left(1+ \sigma \Delta t\right) \rho_j^{n+1} \frac{1}{\sqrt{2\pi}} \sum_{k=1}^{N_v} g_{jk}^{n+1} \omega_k e^{v_k^2/2} = (1+\sigma \Delta t) \rho_j^{n+1}.
\end{align*}
Canceling with $\left(1+ \sigma \Delta t\right) \rho_j^{n+1}$ gives the desired equality for $n+1$, and completes the proof.
\end{proof}

Also, the following inequality is indispensable to show numerical stability of the above system.

\begin{lemma}\label{lemma:fullydiscreteCFL}
Under the time step restriction $\max_k(|v_k|) \Delta t \leq \Delta x$ it holds
\begin{align*}
\Delta t \left\Vert \mathbf{D}^x \mathbf{f}^n \diag\left(v_k\right) - \frac{\Delta x}{2} \mathbf{D}^{xx} \mathbf{f}^n \diag\left(|v_k|\right) \right\Vert_{\mathscr{H}}^2 - \Delta x \left\Vert \mathbf{D}^+ \mathbf{f}^n \diag\left(|v_k|^{1/2}\right)\right\Vert_{\mathscr{H}}^2 \leq 0.
\end{align*}
\end{lemma}
\begin{proof}
We want to apply a Fourier analysis in $x$ that allows us to write the stencil matrices in diagonal form. As in \cite{baumann2024, kusch2023stability}, we define the matrix $\mathbf E
\in\mathbb{C}^{N_x\times N_x}$ as
\begin{align*}
E_{j\alpha} = \sqrt{\frac{\Delta x}{\left| D \right|}}\exp(2\pi i\alpha x_j), \quad j,\alpha = 1,...,N_x,
\end{align*}
where $\left| D \right|$ denotes the length of the domain $D$ and $i\in\mathbb{C}$ the imaginary unit. It is orthonormal, i.e. $\mathbf E\mathbf E^H = \mathbf E^H\mathbf E = \mathbf I$, where the superscript $H$ stands for the complex transpose, and diagonalizes the stencil matrices
\begin{align*}
\mathbf{D}^\gamma\mathbf E = \mathbf E\mathbf{\Lambda}^\gamma, \qquad \text{ with } \gamma \in \{x,xx,+\},
\end{align*}
with $\mathbf{\Lambda}^\gamma
\in\mathbb{C}^{N_x\times N_x}$ being the diagonal matrices with entries
\begin{align*}
\lambda^{x}_{\alpha\alpha} =& \ \frac{1}{2\Delta x}(e^{i\alpha\pi \Delta x}-e^{-i\alpha\pi \Delta x}) = \frac{i}{\Delta x} \sin(\nu_{\alpha}),\\
\lambda_{\alpha\alpha}^{xx} =& \ \frac{1}{\left(\Delta x \right)^2} \left(e^{i\alpha\pi \Delta x}-2+e^{-i\alpha\pi \Delta x}\right) = \frac{2}{\left(\Delta x \right)^2} \left(\cos(\nu_{\alpha})-1\right),\\
\lambda_{\alpha\alpha}^{+} =& \ \frac{1}{\Delta x} \left(e^{i\alpha\pi \Delta x}-1\right) = \frac{1}{\Delta x} \left(\cos(\nu_{\alpha}) + i\sin(\nu_{\alpha})-1\right),
\end{align*}
and $\nu_{\alpha}:=\alpha\pi \Delta x$. Let us denote $\widehat{\mathbf{f}}^n = \left(\widehat{f}_{\alpha k}^n\right) \in \mathbb{C}^{N_x \times N_v}$ with entries $ \widehat{f}_{\alpha k}^n = \sum_{j=1}^{N_x} E_{\alpha j }f_{jk}^n$. With Parseval's identity we obtain
\begin{align*}
& \ \Delta t \left\Vert  \mathbf{D}^x \mathbf{f}^n \diag\left(v_k\right) - \frac{\Delta x}{2} \mathbf{D}^{xx} \mathbf{f}^n \diag\left(|v_k|\right) \right\Vert_{\mathscr{H}}^2  - \Delta x  \left\Vert \mathbf{D}^+ \mathbf{f}^n \diag\left(|v_k|^{1/2}\right)\right\Vert_{\mathscr{H}}^2 \\
=& \ \Delta t \left\Vert  \mathbf{D}^x \mathbf{f}^n \diag\left(v_k \omega_k^{1/2} e^{3v_k^2/4} \right) - \frac{\Delta x}{2} \mathbf{D}^{xx} \mathbf{f}^n \diag\left(|v_k| \omega_k^{1/2} e^{3v_k^2/4}\right)  \right\Vert_F^2\\
&- \Delta x  \left\Vert \mathbf{D}^+ \mathbf{f}^n \diag\left(|v_k|^{1/2} \omega_k^{1/2} e^{3v_k^2/4} \right)\right\Vert_F^2 \\
\overset{\text{Parseval}}{\underset{\text{}}{=}}& \ \Delta t \left\Vert \mathbf{\Lambda}^x \widehat{\mathbf{f}}^n \diag\left(v_k \omega_k^{1/2} e^{3v_k^2/4}\right) - \frac{\Delta x}{2} \mathbf{\Lambda}^{xx} \widehat{\mathbf{f}}^n \diag \left(|v_k| \omega_k^{1/2} e^{3v_k^2/4}) \right) \right\Vert_F^2\\
&- \Delta x \left\Vert \mathbf{\Lambda}^+ \widehat{\mathbf{f}}^n \diag\left(|v_k|^{1/2} \omega_k^{1/2} e^{3v_k^2/4}\right)\right\Vert_F^2\\
=& \ 2 \sum_{\alpha=1}^{N_x} \sum_{k=1}^{N_v} \left( \Delta t \frac{|v_k|^2}{\left(\Delta x\right)^2} (1-\cos(\nu_\alpha)) - \frac{|v_k|}{\Delta x} (1-\cos(\nu_\alpha)) \right) \omega_k e^{3v_k^2/2} \left| \widehat f_{\alpha k}^n \right|^2 .
\end{align*}
A sufficient condition to ensure negativity is that
\begin{align*}
\Delta t \frac{|v_k|^2}{\left(\Delta x\right)^2} (1-\cos(\nu_\alpha)) \leq \frac{|v_k|}{\Delta x} (1-\cos(\nu_\alpha))
\end{align*}
must hold for each index $k$. This leads to the time step restriction $\max_k(|v_k|) \Delta t \leq \Delta x$, which proves the lemma.
\end{proof}

We can now show numerical stability of the proposed system.  

\begin{theorem}\label{theorem:stabilityMg}
Under the time step restriction $\max_k(|v_k|) \Delta t \leq \Delta x$, the fully discrete system \eqref{eqs:KB-both} is numerically stable in the $\mathscr{H}$-norm, i.e. it holds
\begin{align*}
\Vert \mathbf{f}^{n+1}\Vert_{\mathscr{H}}^2 \leq \Vert \mathbf{f}^n\Vert_{\mathscr{H}}^2.
\end{align*}
\end{theorem}

\begin{proof}
We multiply \eqref{eqs:gKB} with  $\rho_j^{n+1} \rho_j^n g_{jk}^{n+1}$ and put the last term of the equation from the right-hand to the left-hand side. This results in
\begin{align*}
\left(\rho_j^{n+1} g_{jk}^{n+1}\right)^2 =& \ \rho_j^n g_{jk}^n \rho_j^{n+1} g_{jk}^{n+1} -\Delta t \sum_{i=1}^{N_x}  \rho_j^{n+1} g_{jk}^{n+1} D_{ji}^x \rho_i^n g_{ik}^n  v_k + \Delta t \frac{\Delta x}{2} \sum_{i=1}^{N_x} \rho_j^{n+1} g_{jk}^{n+1} D_{ji}^{xx} \rho_i^n g_{ik}^n |v_k|\\
&+ \sigma \Delta t \rho_j^{n+1} g_{jk}^{n+1} \left(\rho_j^{n+1} -\rho_j^{n+1} g_{jk}^{n+1}\right).
\end{align*}
Multiplication with $2\left( \frac{1}{\sqrt{2\pi}} e^{-v_k^2/2}\right)^2$ leads to 
\begin{align*}
2\left(f_{jk}^{n+1}\right)^2 =& \ 2 f_{jk}^n f_{jk}^{n+1} -2 \Delta t \sum_{i=1}^{N_x} f_{jk}^{n+1} D_{ji}^x f_{ik}^n  v_k + \Delta t \Delta x \sum_{i=1}^{N_x} f_{jk}^{n+1} D_{ji}^{xx} f_{ik}^n |v_k| + 2 \sigma \Delta t f_{jk}^{n+1} \left(M_{jk}^{n+1} -f_{jk}^{n+1}\right).
\end{align*}
Note that it holds
\begin{align*}
2 f_{jk}^n f_{jk}^{n+1} = \left(f_{jk}^{n+1}\right)^2 +  \left(f_{jk}^n\right)^2 - \left(f_{jk}^{n+1} - f_{jk}^n \right)^2.
\end{align*}
We insert this relation and obtain
\begin{align*}
\left(f_{jk}^{n+1}\right)^2 =& \ \left(f_{jk}^n \right)^2 - \left(f_{jk}^{n+1} - f_{jk}^n \right)^2 - 2\Delta t \sum_{i=1}^{N_x} f_{jk}^{n+1} D_{ji}^x f_{ik}^n  v_k +  \Delta t \Delta x \sum_{i=1}^{N_x} f_{jk}^{n+1} D_{ji}^{xx} f_{ik}^n |v_k|\\
&+ 2 \sigma \Delta t f_{jk}^{n+1} \left(M_{jk}^{n+1} -f_{jk}^{n+1}\right).
\end{align*}
In the next step, we multiply with $\sqrt{2\pi} \omega_k e^{3v_k^2/2}$ and sum over $j$ and $k$. This yields
\begin{align*}
\Vert \mathbf{f}^{n+1} \Vert_{\mathscr{H}}^2 =& \ \Vert \mathbf{f}^n \Vert_{\mathscr{H}}^2 - \sqrt{2\pi} \sum_{j=1}^{N_x} \sum_{k=1}^{N_v}\left(f_{jk}^{n+1} - f_{jk}^n \right)^2 \omega_k e^{3v_k^2/2}\\
&-2 \sqrt{2\pi} \Delta t \sum_{i,j=1}^{N_x} \sum_{k=1}^{N_v} f_{jk}^{n+1} D_{ji}^x f_{ik}^n v_k \omega_k e^{3v_k^2/2} + \sqrt{2\pi} \Delta t \Delta x \sum_{i,j=1}^{N_x} \sum_{k=1}^{N_v} f_{jk}^{n+1} D_{ji}^{xx} f_{ik}^n |v_k| \omega_k e^{3v_k^2/2}\\
&+2 \sqrt{2\pi} \sigma \Delta t \sum_{j=1}^{N_x} \sum_{k=1}^{N_v} f_{jk}^{n+1} \left(M_{jk}^{n+1} -f_{jk}^{n+1}\right) \omega_k e^{3v_k^2/2}.
\end{align*}
According to Lemma \ref{lemma:fullydiscrete-rho}, we can use the equality $\sum_{k=1}^{N_v} f_{jk}^{n+1} \omega_k e^{v_k^2} = \rho_j^{n+1}$. Hence, we can conclude that the term $2 \sqrt{2\pi} \sigma \Delta t \sum_{j=1}^{N_x} \sum_{k=1}^{N_v} M_{jk}^{n+1} \left(M_{jk}^{n+1} - f_{jk}^{n+1} \right) \omega_k e^{3v_k^2/2}$ is equal to zero. Lemma \ref{lemma:stencilmatrices} implies that also the term $2 \sqrt{2\pi} \Delta t \sum_{i,j=1}^{N_x} \sum_{k=1}^{N_v} f_{jk}^n D_{ji}^x f_{ik}^n v_k \omega_k e^{3v_k^2/2}$ is equal to zero. We subtract both and add an additional zero by adding and subtracting the second-order term $\sqrt{2\pi} \Delta t \Delta x \sum_{i,j=1}^{N_x} \sum_{k=1}^{N_v} f_{jk}^n D_{ji}^{xx} f_{ik}^n |v_k| \omega_k e^{3v_k^2/2}$. This leads to
\begin{align*}
\Vert \mathbf{f}^{n+1} \Vert_{\mathscr{H}}^2 =& \ \Vert \mathbf{f}^n \Vert_{\mathscr{H}}^2 - \sqrt{2\pi} \sum_{j=1}^{N_x} \sum_{k=1}^{N_v}\left(f_{jk}^{n+1} - f_{jk}^n \right)^2 \omega_k e^{3v_k^2/2}\\
&-2 \sqrt{2\pi} \Delta t \sum_{i,j=1}^{N_x} \sum_{k=1}^{N_v} \left(f_{jk}^{n+1} - f_{jk}^n \right) D_{ji}^x f_{ik}^n v_k \omega_k e^{3v_k^2/2}\tag{I}\\
&+ \sqrt{2\pi} \Delta t \Delta x \sum_{i,j=1}^{N_x} \sum_{k=1}^{N_v} \left( f_{jk}^{n+1} - f_{jk}^n \right) D_{ji}^{xx} f_{ik}^n |v_k| \omega_k e^{3v_k^2/2}\tag{II}\\
&+ \sqrt{2\pi} \Delta t \Delta x \sum_{i,j=1}^{N_x} \sum_{k=1}^{N_v} f_{jk}^n D_{ji}^{xx} f_{ik}^n |v_k| \omega_k e^{3v_k^2/2}\tag{III}\\ 
&- 2 \sqrt{2\pi} \sigma \Delta t \sum_{j=1}^{N_x} \sum_{k=1}^{N_v} \left( f_{jk}^{n+1} - M_{jk}^{n+1} \right)^2 \omega_k e^{3v_k^2/2}.
\end{align*}
Now, we analyze the terms (I), (II) and (III) separately. Let us start with (I) and (II) and apply Young's inequality which states that for $a,b \in \mathbb{R}$ it holds that $a \cdot b \leq \frac{a^2}{2} + \frac{b^2}{2}$. For the sum (I) +(II) this yields
\begin{align*}
&-2 \sqrt{2\pi} \Delta t \sum_{i,j=1}^{N_x} \sum_{k=1}^{N_v} \left(f_{jk}^{n+1} - f_{jk}^n \right) D_{ji}^x f_{ik}^n v_k \omega_k e^{3v_k^2/2} + \sqrt{2\pi} \Delta t \Delta x \sum_{i,j=1}^{N_x} \sum_{k=1}^{N_v} \left( f_{jk}^{n+1} - f_{jk}^n \right) D_{ji}^{xx} f_{ik}^n |v_k| \omega_k e^{3v_k^2/2}\\
=& \sum_{j=1}^{N_x} \sum_{k=1}^{N_v} \left(-\sqrt{2} \sqrt[4]{2\pi} \left(f_{jk}^{n+1} - f_{jk}^n \right) \sqrt{\omega_k} e^{3v_k^2/4}\right) \left( \sqrt{2} \sqrt[4]{2\pi} \Delta t \sum_{i=1}^{N_x} \left( D_{ji}^x f_{ik}^n v_k - \frac{\Delta x}{2} D_{ji}^{xx} f_{ik}^n |v_k| \right) \sqrt{\omega_k} e^{3v_k^2/4} \right)\\
\leq& \ \sqrt{2\pi}\sum_{j=1}^{N_x} \sum_{k=1}^{N_v} \left(f_{jk}^{n+1} - f_{jk}^n \right)^2 \omega_k e^{3v_k^2/2} + \sqrt{2\pi} \left(\Delta t\right)^2 \sum_{j=1}^{N_x} \sum_{k=1}^{N_v} \left(\sum_{i=1}^{N_x} \left( D_{ji}^x f_{ik}^n v_k - \frac{\Delta x}{2} D_{ji}^{xx} f_{ik}^n |v_k| \right) \right)^2 \omega_k e^{3v_k^2/2}.
\end{align*}
For (III) we exploit the properties of the spatial stencil matrices given in Lemma \ref{lemma:stencilmatrices}. We derive the equality 
\begin{align*}
\sqrt{2\pi} \Delta t \Delta x \sum_{i,j=1}^{N_x} \sum_{k=1}^{N_v} f_{jk}^n D_{ji}^{xx} f_{ik}^n |v_k| \omega_k e^{3v_k^2/2} = - \sqrt{2\pi} \Delta t \Delta x \sum_{j=1}^{N_x} \sum_{k=1}^{N_v} \left(\sum_{i=1}^{N_x} D_{ji}^+ f_{ik}^n |v_k|^{1/2} \right)^2\omega_k e^{3v_k^2/2}. 
\end{align*}
We insert both relations and obtain
\begin{align*}
\Vert \mathbf{f}^{n+1} \Vert_{\mathscr{H}}^2 \leq& \ \Vert \mathbf{f}^n \Vert_{\mathscr{H}}^2 + \sqrt{2\pi} \left(\Delta t\right)^2 \sum_{i=1}^{N_x} \sum_{k=1}^{N_v} \left(\sum_{j=1}^{N_x} D_{ji}^x f_{ik}^n v_k - \frac{\Delta x}{2} D_{ji}^{xx} f_{ik}^n |v_k| \right)^2 \omega_k e^{3v_k^2/2}\\
&- \sqrt{2\pi} \Delta t \Delta x  \sum_{j=1}^{N_x} \sum_{k=1}^{N_v} \left(\sum_{i=1}^{N_x} D_{ji}^+ f_{ik}^n |v_k|^{1/2} \right)^2\omega_k e^{3v_k^2/2}\\ 
&- 2 \sqrt{2\pi} \sigma \Delta t \sum_{j=1}^{N_x} \sum_{k=1}^{N_v} \left( f_{jk}^{n+1} - M_{jk}^{n+1} \right)^2 \omega_k e^{3v_k^2/2}.
\end{align*}
With Lemma \ref{lemma:fullydiscreteCFL} we can conclude that under the CFL condition $\max_k(|v_k|)\Delta t \leq \Delta x$ it holds $\Vert \mathbf{f}^{n+1} \Vert_{\mathscr{H}}^2 \leq \ \Vert \mathbf{f}^n \Vert_{\mathscr{H}}^2$. Hence, under this time step restriction the proposed fully discrete system \eqref{eqs:KB-both} is numerically stable in the $\mathscr{H}$-norm.
\end{proof}

\section{Dynamical low-rank approximation for the stable $Mg$ system}\label{sec5:DLRAforMg}

In practical applications, the implementation of the full system given in \eqref{eqs:KB-both} may lead to prohibitive numerical costs, especially when computing in higher-dimensional settings. To reduce computational and memory demands, we apply dynamical low-rank approximation to the distribution function $g$.

\subsection{Background on dynamical low-rank approximation}\label{sec5-1:generalDLRA}

The concept of dynamical low-rank approximation has been introduced in a semi-discrete time-dependent matrix setting \cite{koch2007dynamical}. Let us consider $\mathbf{g} \in \mathbb{R}^{N_x \times N_v}$ being the solution of the matrix differential equation
\begin{align*}
\dot{\mathbf{g}}(t) = \mathbf{F}\left(\mathbf{g}(t) \right),
\end{align*}
where $\mathbf{F}: \mathbb{R}^{N_x \times N_v} \to \mathbb{R}^{N_x \times N_v}$ denotes the right-hand side of the equation. We then seek for an approximation of $\mathbf{g}$ in the following form
\begin{align}\label{sec5:DLRAmatrixform}
\mathbf{g}_r(t) = \mathbf{X}(t) \mathbf{S}(t) \mathbf{V}(t)^\top,
\end{align}
with $\mathbf{X} \in \mathbb{R}^{N_x \times r}$ and $\mathbf{V} \in \mathbb{R}^{N_v \times r}$ denoting the orthonormal spatial and orthonormal velocity basis, respectively. The slim matrix $\mathbf{S} \in \mathbb{R}^{r \times r}$ is called the coefficient or coupling matrix and determines the rank $r$ of the approximation. The set of all matrices of the above form \eqref{sec5:DLRAmatrixform} constitute the low-rank manifold $\mathcal{M}_r$. Its corresponding tangent space at $\mathbf{g}_r(t)$ shall be denoted by $\mathcal{T}_{\mathbf{g}_r(t)}\mathcal{M}_r$. We now look for $\mathbf{g}_r(t) \in \mathcal{M}_r$ such that at all times $t$ the minimization problem 
\begin{align*}
\min_{\dot{\mathbf{g}}_r(t) \in \mathcal{T}_{\mathbf{g}_r(t)} \mathcal{M}_r} \Vert \dot{\mathbf{g}}_r(t) - \mathbf{F}\left(\mathbf{g}_r(t)\right)\Vert_F
\end{align*}
is fulfilled. Following \cite{koch2007dynamical}, this minimization constraint is equivalent to determining $\dot{\mathbf{g}}_r(t) \in \mathcal{T}_{\mathbf{g}_r(t)}\mathcal{M}_r$ by an orthogonal projection onto the tangent space such that
\begin{align}\label{sec5:projection}
\dot{\mathbf{g}}_r(t) = \mathbf{P}(\mathbf{g}_r(t)) \mathbf{F}(\mathbf{g}_r(t)),
\end{align}
where the orthogonal projector $\mathbf{P}$ onto $\mathcal{T}_{\mathbf{g}_r(t)}\mathcal{M}_r$ applied to an arbitrary quantity $\mathbf{G}$ can explicitly be given as 
\begin{align*}
\mathbf{P}(\mathbf{g}_r(t)) \mathbf{G} = \mathbf{X} \mathbf{X}^\top \mathbf{G} - \mathbf{X} \mathbf{X}^\top \mathbf{G} \mathbf{V} \mathbf{V}^\top + \mathbf{G} \mathbf{V} \mathbf{V}^\top.
\end{align*}
Different robust time integrators to solve \eqref{sec5:projection} exist. Frequently used integrators are the projector-splitting integrator \cite{lubich2014projector}, the basis update \& Galerkin integrator \cite{ceruti2022unconventional} as well as its rank-adaptive extension \cite{ceruti2022rank} and the parallel integrator \cite{ceruti2024parallel}. In this paper, we make use of the rank-adaptive BUG integrator whose concept shall be explained in the following. 

In the first two steps, the rank-adaptive BUG integrator updates and augments the bases $\mathbf{X}$ and $\mathbf{V}$ in parallel such that their rank increases from $r$ to $2r$, for the spatial and velocity basis respectively. We denote the augmented quantities of rank $2r$ with hats. Then, for the augmented bases a Galerkin step is performed before in a last step a new rank $r_1 \leq 2r$ is determined using a truncation with prescribed tolerance. In detail, the rank-adaptive BUG integrator performs the following steps in order to update the matrix $\mathbf{g}_r^n = \mathbf{X}^n \mathbf{S}^n \mathbf{V}^{n,\top}$ at time $t_n$ to $\mathbf{g}_r^{n+1} = \mathbf{X}^{n+1} \mathbf{S}^{n+1} \mathbf{V}^{n+1,\top}$ at time $t_{n+1}=t_n + \Delta t$: 

\textbf{\textit{K}-Step}: Let us fix the velocity basis $\mathbf{V}^n$ at time $t_n$ and introduce the notation $\mathbf{K}(t) = \mathbf{X}(t) \mathbf{S}(t)$. The spatial basis $\mathbf{X}^n$ is updated and augmented by first solving the PDE
\begin{align*}
\dot{\mathbf{K}}(t) = \mathbf{F}\left( \mathbf{K}(t) \mathbf{V}^{n,\top} \right) \mathbf{V}^n, \quad \mathbf{K}(t_n) = \mathbf{X}^n \mathbf{S}^n,
\end{align*}
and then determining $\widehat{\mathbf{X}}^{n+1} \in \mathbb{R}^{N_x \times 2r}$ as an orthonormal basis of $[\mathbf{K}(t_{n+1}), \mathbf{X}^n] \in \mathbb{R}^{N_x \times 2r}$, e.g. by QR-decomposition. Then, we compute $\widehat{\mathbf{M}} = \widehat{\mathbf{X}}^{n+1,\top} \mathbf{X}^n \in \mathbb{R}^{2r \times r}$.

\textbf{\textit{L}-Step}:
Let us fix the spatial basis $\mathbf{X}^n$ at time $t_n$ and introduce the notation $\mathbf{L}(t) = \mathbf{V}(t) \mathbf{S}(t)^\top$. The velocity basis $\mathbf{V}^n$ is updated and augmented by first solving the PDE
\begin{align*}
\dot{\mathbf{L}}(t) = \mathbf{F}\left(\mathbf{X}^n \mathbf{L}(t)^\top \right)^\top \mathbf{X}^n, \quad \mathbf{L}(t_n) = \mathbf{V}^n \mathbf{S}^{n,\top},
\end{align*}
and then determining $\widehat{\mathbf{V}}^{n+1} \in \mathbb{R}^{N_v \times 2r}$ as an orthonormal basis of $[\mathbf{L}(t_{n+1}), \mathbf{V}^n] \in \mathbb{R}^{N_v \times 2r}$, e.g. by QR-decomposition. Then, we compute $\widehat{\mathbf{N}} = \widehat{\mathbf{V}}^{n+1,\top} \mathbf{V}^n \in \mathbb{R}^{2r \times r}$.

\textbf{\textit{S}-step}: Update the coupling matrix from $\mathbf{S}^n \in \mathbf{R}^{r \times r}$ to $\widehat{\mathbf{S}}^{n+1} \in \mathbb{R}^{2r \times 2r}$ by solving the ODE
\begin{align*}
\dot{\widehat{\mathbf{S}}}(t) = \widehat{\mathbf{X}}^{n+1, \top} \mathbf{F} \left( \widehat{\mathbf{X}}^{n+1} \widehat{\mathbf{S}}(t) \widehat{\mathbf{V}}^{n+1,\top} \right) \widehat{\mathbf{V}}^{n+1}, \quad \widehat{\mathbf{S}}(t_n) = \widehat{\mathbf{M}} \mathbf{S}^n \widehat{\mathbf{N}}^\top.
\end{align*}

\textbf{Truncation}: Compute the singular value decomposition of $\widehat{\mathbf{S}}^{n+1} = \widehat{\mathbf{P}} \mathbf{\Sigma}\widehat {\mathbf{Q}}^\top$ with $\mathbf{\Sigma} = \text{diag}(\sigma_j)$. For a prescribed tolerance parameter $\vartheta$ the new rank $r_1 \leq 2r$ is chosen such that 
\begin{align*}
\left(\sum_{j=r_1+1}^{2r} \sigma_j^2\right)^{1/2} \leq \vartheta.
\end{align*}
Let now $\mathbf{S}^{n+1} \in \mathbb{R}^{r_1\times r_1}$ contain the $r_1$ largest singular values and $\mathbf{P}^{n+1} \in \mathbb{R}^{2r\times r_1}$ and $\mathbf{Q}^{n+1} \in \mathbb{R}^{2r\times r_1}$ contain the first $r_1$ columns of $\widehat{\mathbf{P}}$ and $\widehat {\mathbf{Q}}$, respectively. Then, the time-updated spatial basis can be determined as $\mathbf{X}^{n+1} = \widehat{\mathbf{X}}^{n+1}\mathbf{P}^{n+1} \in \mathbb{R}^{N_x \times r_1}$ and the time-updated velocity basis as $\mathbf{V}^{n+1} = \widehat{\mathbf{V}}^{n+1} \mathbf{Q}^{n+1} \in \mathbb{R}^{N_v \times r_1}$.

Altogether, the time-updated approximation of the distribution function after one time step is given by $\mathbf{g}_r^{n+1} = \mathbf{X}^{n+1} \mathbf{S}^{n+1} \mathbf{V}^{n+1,\top}$.

\subsection{DLRA algorithm for multiplicative linear Boltzmann-BGK}

In this section, we apply a DLRA approach to the stable and fully discretized system \eqref{eqs:KB-both}. Putting all terms containing $g_{jk}^{n+1}$ to the left-hand side of \eqref{eqs:gKB} and multiplying it with $\frac{\rho_j^n}{\rho_j^{n+1}}$, equations \eqref{eqs:KB-both} can be written in the equivalent form
\begin{subequations}
\begin{align}
g_{jk}^{n+1} \left(1+\sigma \Delta t \right) =& \ \frac{\rho_j^n}{\rho_j^{n+1}} g_{jk}^n -\Delta t \sum_{i=1}^{N_x} \frac{1}{\rho_j^{n+1}} D_{ji}^x \left(\rho_i^n g_{ik}^n \right) v_k + \Delta t \frac{\Delta x}{2} \sum_{i=1}^{N_x} \frac{1}{\rho_j^{n+1}} D_{ji}^{xx} \left(\rho_i^n g_{ik}^n \right) |v_k| + \sigma \Delta t,\label{eqs:KBrefDLRA-g}\\
\rho_j^{n+1} =& \ \rho_j^n - \Delta t \frac{1}{\sqrt{2\pi}} \sum_{i=1}^{N_x} \sum_{k=1}^{N_v} D_{ji}^x \rho_i^n g_{ik}^n v_k \omega_k e^{v_k^2/2} + \Delta t \frac{\Delta x}{2\sqrt{2\pi}} \sum_{i=1}^{N_x} \sum_{k=1}^{N_v} D_{ji}^{xx} \rho_i^n g_{ik}^n |v_k| \omega_k e^{v_k^2/2}.\label{eqs:KBrefDLRA-rho}
\end{align}\label{eqs:KBrefDLRA-both}
\end{subequations}

We propose a numerically stable DLRA implementation that uses the rank-adaptive BUG integrator presented in \cite{ceruti2022rank} and introduced in the previous subsection for equation \eqref{eqs:KBrefDLRA-g} together with additional basis augmentations and a suitable truncation strategy. Note that in the following, for simplicity, we write $\mathbf{g} = \left( g_{jk} \right)$ instead of $\mathbf{g}_r$.

Starting from \eqref{eqs:KBrefDLRA-both}, we apply the rank-adaptive BUG integrator, leading to a splitting of \eqref{eqs:KBrefDLRA-g} into three substeps, the $K$-, $L$-, and $S$-step. In the $K$- as well as in the $L$-step we perform an additional basis augmentation ensuring certain quantities to be contained in the basis at all times. The explicit form of the basis augmentations will be made clear later in the proof of stability of the proposed scheme. The augmented bases are used to determine the $S$-step. Note that the scattering term $\left(1+\sigma \Delta t \right)$ is only applied in the $S$-step as it does not affect the span of the basis functions derived in the $K$- and $L$-step. We obtain the following DLRA scheme:

\begin{subequations}
We first substitute $g_{jk}^n = \sum_{m,\ell=1}^r X_{jm}^n S_{m\ell}^n V_{k\ell}^n$ into the update equation \eqref{eqs:KBrefDLRA-rho} and obtain
\begin{align}
\rho_j^{n+1} =& \ \rho_j^n - \Delta t \frac{1}{\sqrt{2\pi}} \sum_{i=1}^{N_x} D_{ji}^x \rho_i^n \sum_{m,\ell=1}^r X_{im}^n S_{m\ell}^n \sum_{k=1}^{N_v}  V_{k\ell}^n v_k \omega_k e^{v_k^2/2}\label{alg:rho}\\
&+ \Delta t \frac{\Delta x}{2\sqrt{2\pi}} \sum_{i=1}^{N_x} D_{ji}^{xx} \rho_i^n \sum_{m,\ell=1}^r X_{im}^n S_{m\ell}^n \sum_{k=1}^{N_v} V_{k\ell}^n |v_k| \omega_k e^{v_k^2/2}.\nonumber
\end{align}
For the $K$-step we introduce the notation $K_{j\ell}^n = \sum_{m=1}^r X_{jm}^n S_{m\ell}^n$ and solve
\begin{align}
K_{jp}^{n+1} =& \ \frac{\rho_j^n}{\rho_j^{n+1}} K_{jp}^n -\Delta t \frac{1}{\rho_j^{n+1}} \sum_{i=1}^{N_x} D_{ji}^x \rho_i^n \sum_{\ell=1}^r K_{i\ell}^n \sum_{k=1}^{N_v} V_{k\ell}^n v_k V_{kp}^n\label{alg:Kstep}\\
&+ \Delta t \frac{\Delta x}{2} \frac{1}{\rho_j^{n+1}} \sum_{i=1}^{N_x} D_{ji}^{xx} \rho_i^n \sum_{\ell=1}^r K_{i\ell}^n \sum_{k=1}^{N_v} V_{k\ell}^n |v_k| V_{kp}^n+ \sigma \Delta t \sum_{k=1}^{N_v} V_{kp}^n.\nonumber
\end{align}
This gives the updated matrix $ \mathbf{K}^{n+1} = (K_{jp}^{n+1})$ to which together with the old basis $\mathbf{X}^n = (X_{jm}^n)$ a QR-decomposition is applied, giving the new augmented basis $\widehat{\mathbf{X}}^{n+1} = (\widehat{X}_{jm}^{n+1})$.

In addition, we augment this basis according to
\begin{align}
\doublehat{\mathbf{X}}^{n+1}  = \text{qr} \left( \left[\widehat{\mathbf{X}}^{n+1}, \left(\boldsymbol{\rho}^{n+1}\right)^2 \widehat{\mathbf{X}}^{n+1} \right] \right)\label{alg:KbasisAug}
\end{align}
leading to a new augmented basis $\doublehat{\mathbf{X}}^{n+1} = (\doublehat{X}_{jm}^{n+1})$ of rank $4r$, and we compute $\doublehat{\mathbf{M}} = \doublehat{\mathbf{X}}^{n+1,\top} \mathbf{X}^n$. Note that all quantities of rank $2r$ are denoted with one hat and all quantities of rank $4r$ with double hats. 

For the $L$-step we write $L_{mk}^n = \sum_{\ell=1}^r S_{\ell m}^n V_{\ell k}^n$ and solve
\begin{align}
L_{pk}^{n+1}=& \ \sum_{m=1}^r L_{mk}^n \sum_{j=1}^{N_x} X_{mj}^n \frac{\rho_j^n}{\rho_j^{n+1}} X_{pj}^n - \Delta t \sum_{m=1}^r v_k L_{mk}^n \sum_{i=1}^{N_x} X_{mi}^n \rho_i^n \sum_{j=1}^{N_x} D_{ij}^x \frac{1}{\rho_j^{n+1}} X_{pj}^n\label{alg:Lstep}\\
&+ \Delta t \frac{\Delta x}{2} \sum_{m=1}^r |v_k| L_{mk}^n \sum_{i=1}^{N_x} X_{mi}^n \rho_i^n \sum_{j=1}^{N_x} D_{ij}^{xx} \frac{1}{\rho_j^{n+1}} X_{pj}^n + \sigma \Delta t \sum_{j=1}^{N_x} X_{pj}^n.\nonumber
\end{align}
This gives the updated matrix $ \mathbf{L}^{n+1} = (L_{pk}^{n+1})$ to which together with the old basis $\mathbf{V}^n = (V_{\ell k}^n)$ a QR-decomposition is applied, giving the new augmented basis $\widehat{\mathbf{V}}^{n+1} = (\widehat{V}_{\ell k}^{n+1})$. 

In addition, we augment this basis according to
\begin{align}
\doublehat{\mathbf{V}}^{n+1}  = \text{qr} \left(\left[\widehat{\mathbf{V}}^{n+1}, \mathbf{\omega} e^{\mathbf{v}^2/2} \widehat{\mathbf{V}}^{n+1} \right]\right) \label{alg:LbasisAug}
\end{align}
leading to a new augmented basis $\doublehat{\mathbf{V}}^{n+1} = (\doublehat{V}_{\ell k}^{n+1})$ of rank $4r$, and we compute $\doublehat{\mathbf{N}} = \doublehat{\mathbf{V}}^{n+1,\top} \mathbf{V}^n$.

For the $S$-step we denote $\widetilde S_{m\ell}^n = \sum_{j,k=1}^r \doublehat M_{mj} S_{jk}^n \doublehat N_{\ell k}$ and insert the expressions $g_{jk}^n = \sum_{m,\ell=1}^{4r} \doublehat X_{jm}^{n+1} \widetilde S_{m\ell}^n \doublehat V_{k\ell}^{n+1}$ and $g_{jk}^{n+1} = \sum_{m,\ell=1}^{4r} \doublehat X_{jm}^{n+1} \doublehat S_{m\ell}^{n+1} \doublehat V_{k\ell}^{n+1}$ into \eqref{eqs:KBrefDLRA-g}. We multiply with $\doublehat X_{jq}^{n+1} \doublehat V_{kp}^{n+1}$ and sum over $j$ and $k$. This leads to 
\begin{align}
\doublehat S_{qp}^{n+1} =\ \frac{1}{ 1+\sigma \Delta t} \Bigg(& \sum_{j=1}^{N_x} \doublehat X_{jq}^{n+1} \frac{\rho_j^n}{\rho_j^{n+1}} \sum_{m,\ell=1}^{4r} \doublehat X_{jm}^{n+1} \widetilde S_{m\ell}^n \sum_{k=1}^{N_v} \doublehat V_{k\ell}^{n+1} \doublehat V_{kp}^{n+1}\nonumber\\
&- \Delta t \sum_{j=1}^{N_x} \doublehat X_{jq}^{n+1} \frac{1}{\rho_j^{n+1}} \sum_{i=1}^{N_x} D_{ji}^x \rho_i^n \sum_{m,\ell=1}^{4r} \doublehat X_{im}^{n+1} \widetilde S_{m\ell}^n \sum_{k=1}^{N_v} \doublehat V_{k\ell}^{n+1} v_k \doublehat V_{kp}^{n+1}\label{alg:Sstep}\\
&+ \Delta t \frac{\Delta x}{2} \sum_{j=1}^{N_x} \doublehat X_{jq}^{n+1} \frac{1}{\rho_j^{n+1}} \sum_{i=1}^{N_x} D_{ji}^{xx} \rho_i^n \sum_{m,\ell=1}^{4r} \doublehat X_{im}^{n+1} \widetilde S_{m\ell}^n \sum_{k=1}^{N_v}  \doublehat V_{k\ell}^{n+1} |v_k| \doublehat V_{kp}^{n+1}\nonumber\\
&+ \sigma \Delta t  \sum_{j=1}^{N_x} \doublehat X_{jq}^{n+1} \sum_{k=1}^{N_v} \doublehat V_{kp}^{n+1}\Bigg).\nonumber
\end{align}

The last step consists in truncating the augmented quantities $\doublehat{\mathbf{X}}^{n+1}$, $\doublehat{\mathbf{V}}^{n+1}$ and $\doublehat{\mathbf{S}}^{n+1}$ from rank $4r$ to a new rank $r_1$. We use a modification of the truncation strategy described in Section \ref{sec5-1:generalDLRA}
that ensures that the equality $\frac{1}{\sqrt{2\pi}} \sum_{\ell,m=1}^r \sum_{k=1}^{N_v} X_{jm}^n S_{m\ell}^n V_{k\ell}^n \omega_k e^{v_k^2/2} = 1$ stays valid in each time step and works as follows: 

Let us denote $\mathbf{Z} \in \mathbb{R}^{N_v}$ being the vector with entries $Z_k = \frac{1}{\sqrt{2 \pi}}\omega_k e^{v_k^2/2}$ and let $\mathbf{z} = \frac{\mathbf{Z}}{\Vert \mathbf{Z} \Vert_E}$, where $\Vert \cdot \Vert_E$ stands for the Euclidean norm. We then want to have
\begin{align*}
\mathbf{1} = \doublehat{\mathbf{X}}^{n+1} \doublehat{\mathbf{S}}^{n+1} \doublehat{\mathbf{V}}^{n+1,\top} \mathbf{Z} = \left( \doublehat{\mathbf{X}}^{n+1} \doublehat{\mathbf{S}}^{n+1} \doublehat{\mathbf{V}}^{n+1,\top} \mathbf{z} \mathbf{z}^\top + \doublehat{\mathbf{X}}^{n+1} \doublehat{\mathbf{S}}^{n+1} \doublehat{\mathbf{V}}^{n+1,\top} \left( \mathbf{I} - \mathbf{z} \mathbf{z}^\top \right)\right) \mathbf{Z} = : \left(\mathbf{H}_1 + \mathbf{H}_2 \right) \mathbf{Z},
\end{align*}
with $\mathbf{H}_1 = \doublehat{\mathbf{X}}^{n+1} \doublehat{\mathbf{S}}^{n+1} \doublehat{\mathbf{V}}^{n+1,\top} \mathbf{z} \mathbf{z}^\top, \mathbf{H}_2 = \doublehat{\mathbf{X}}^{n+1} \doublehat{\mathbf{S}}^{n+1} \doublehat{\mathbf{V}}^{n+1,\top} \left( \mathbf{I} - \mathbf{z} \mathbf{z}^\top \right)$, $\mathbf{I} \in \mathbb{R}^{N_v \times N_v}$ denoting the identity matrix and $\mathbf{1} \in \mathbb{R}^{N_x}$ the vector containing ones at each entry. $\mathbf{H}_1$ is a matrix of rank $1$. We determine its low-rank factors by a singular value decomposition such that $\mathbf{X}^{\mathbf{H}_1} \mathbf{S}^{\mathbf{H}_1} \mathbf{V}^{\mathbf{H}_1,\top}= \text{svd}\left(\doublehat{\mathbf{S}}^{n+1} \doublehat{\mathbf{V}}^{n+1,\top} \mathbf{z} \mathbf{z}^\top \right)$ with $\mathbf{X}^{\mathbf{H}_1} \in \mathbb{R}^{4r}, \mathbf{S}^{\mathbf{H}_1} \in \mathbb{R}$, and $\mathbf{V}^{\mathbf{H}_1} \in \mathbb{R}^{N_v}$. For $\mathbf{H}_2$, it holds that $\mathbf{H}_2 \mathbf{Z} =0$. We apply the truncation strategy from Section \ref{sec5-1:generalDLRA} to $\mathbf{H}_2$ and obtain $\mathbf{X}^\ast, \mathbf{S}^\ast,$ and $\mathbf{V}^\ast$, which shall be of rank $\widetilde r_1$. Finally, we combine both parts by performing a QR-decomposition of
\begin{align*}
\mathbf{X}^{n+1} \mathbf{R}^1 = \left[\doublehat{\mathbf{X}}^{n+1} \mathbf{X}^{\mathbf{H}_1}, \mathbf{X}^\ast \right], \quad \text{ and } \quad \mathbf{V}^{n+1} \mathbf{R}^2 = \left[ \mathbf{V}^{\mathbf{H}_1}, \mathbf{V}^\ast \right],
\end{align*}
and setting
\begin{align*}
\mathbf{S}^{n+1} = \mathbf{R}^1 \begin{bmatrix}
\mathbf{S}^{\mathbf{H}_1} & 0 \\ 
0 & \mathbf{S}^\ast
\end{bmatrix} \mathbf{R}^{2,\top}.
\end{align*}
The new rank $r_1$ is then given by $r_1 = \widetilde r_1 + 1$. With the time-updated low-rank factors $\mathbf{X}^{n+1}$, $\mathbf{V}^{n+1}$ and $\mathbf{S}^{n+1}$, the updated low-rank approximation of the solution is $g_{jk}^{n+1} = \sum_{m,\ell=1}^{r_1} X_{jm}^{n+1} S_{m\ell}^{n+1} V_{k\ell}^{n+1}$. The steps of the proposed DLRA scheme are visualized in Figure \ref{fig:flowchart}. Note that the notation using brackets refers to a simplification of the algorithm that is explained later in Section \ref{sec:1DPlanesource}.\label{alg:DLRA}
\end{subequations}

\begin{figure}[h!]
    \centering
    \begin{tikzpicture}[scale=0.95, node distance = 3.5cm,auto,font=\sffamily]
        \node[block](init){
        \textbf{input}
    \begin{varwidth}{\linewidth}\begin{itemize}
        \item density at time $t_n$: $\rho_j^n$
        \item low-rank factors at time $t_n$:  $X^n_{jm},S^n_{m\ell}, V^n_{k\ell}$
        \item rank at time $t_n$: $r$
    \end{itemize}\end{varwidth}
        };
        \node[block, below = 0.35cm of init, fill=blue!20](rhostep){update density according to \eqref{alg:rho}};
        \node[block, below = 0.75cm of rhostep, fill=blue!20](KLstep){update bases according to \eqref{alg:Kstep} and \eqref{alg:Lstep}};
        \node[block, below = 0.75cm of KLstep,fill=red!20](augment1){augment bases with $X^n_{jm}, V^n_{k\ell}$};
        \node[block, below = 0.75cm of augment1,fill=red!20](augment2){augment bases with $\left(\rho_j^{n+1}\right)^2 \widehat X^{n+1}_{jm}$ and $ \omega_k e^{v_k^2/2} \widehat V^{n+1}_{k\ell}$ according to \eqref{alg:KbasisAug} and \eqref{alg:LbasisAug}};
        \node[block, below = 0.75cm of augment2,fill=blue!20](SStep){update coefficient matrix according to \eqref{alg:Sstep}};
        \node[block, below = 0.75cm of SStep,fill=yellow!20](truncate){truncate factors $\doublehat {X}^{n+1}_{jm},\doublehat{S}^{n+1}_{m\ell},\doublehat{V}^{n+1}_{k\ell}$ (or $\widehat {X}^{n+1}_{jm},\widehat{S}^{n+1}_{m\ell},\widehat{V}^{n+1}_{k\ell}$ )};
        \node[block, below = 0.35cm of truncate](out){
        \textbf{output}
    \begin{varwidth}{\linewidth}\begin{itemize}
        \item density at time $t_{n+1}$: $\rho_j^{n+1}$
        \item low-rank factors at time $t_{n+1}$:  $X^{n+1}_{jm},S^{n+1}_{m\ell},V^{n+1}_{k\ell}$
        \item rank at time $t_{n+1}$: $r_1$
    \end{itemize}\end{varwidth}
        };
        \path[line] (init) -- node [near end] {} (rhostep);
        \path[line] (rhostep) -- node [near end] {} (KLstep);
        \path[line] (KLstep) -- node [near end] {} (augment1);
        \path[line] (augment1) -- node [near end] {} (augment2);
        \path[line] (augment2) -- node [near end] {} (SStep);
        \path[line] (init.east) -- ([xshift=0.35cm] init.east) |- (augment1.east);
        \path[line] (augment1.west) -- ([xshift=-4.4cm] augment1.west) |- (SStep.west);
        \path[line] (SStep) -- node [near end] {} (truncate);
        \path[line] (truncate) -- node [near start] {} (out);
        \node[below right=0.0cm and 0.1cm of augment1.south](test){$\widehat{X}_{jm}^{n+1}, \widehat{V}_{k\ell}^{n+1}$};
        \node[above=0.1cm of augment2, xshift=-5.2cm](test){($\widehat{X}_{jm}^{n+1}, \widehat{V}_{k\ell}^{n+1}$)};
        \node[below right=0.0cm and 0.1cm of augment2.south](test){$\doublehat{X}_{jm}^{n+1}, \doublehat{V}_{k\ell}^{n+1}$};
        \node[below right=0.0cm and 0.1cm of SStep.south](test){$\doublehat S^{n+1}_{m\ell}$ (or $\widehat S^{n+1}_{m\ell}$)};
        \node[below right=0.0cm and 0.1cm of rhostep.south](test){$\rho_j^{n+1}$};
        \node[below right=0.0cm and 0.1cm of KLstep.south](test){$K^{n+1}_{jm}, L^{n+1}_{k\ell}$};
    \end{tikzpicture}
    \caption{Flowchart of the (simplified) stable DLRA scheme \eqref{alg:DLRA}.}
    \label{fig:flowchart}
\end{figure}
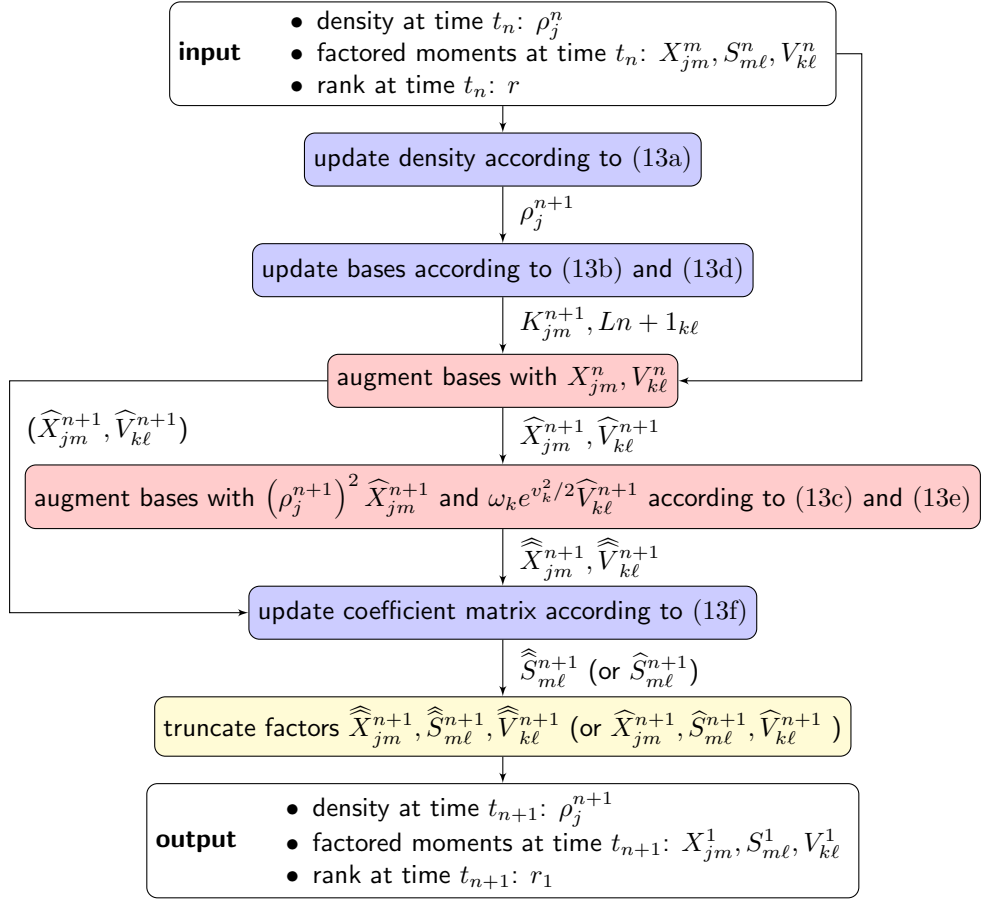

% We have to ensure that $\rho = \int f \mathrm{d}v$ holds over all time steps. Hence, the truncation has to be adjusted. 

% We need $\frac{1}{\sqrt{2\pi}} \sum_{k=1}^{N_v} g_{jk} \omega_k e^{v_k^2/2} = 1$. This can be equivalently rewritten as $\mathbf{g} \mathbf{W} = 1$ with $\mathbf{W} = \frac{1}{\sqrt{2\pi}} \begin{pmatrix}
%     \omega_1 e^{v_1^2/2} \\ \omega_2 e^{v_2^2/2} \\ ...
% \end{pmatrix}$.
% Let us denote the normalized vector $\mathbf{W}$ as $ \mathbf{e}$. Then we want to have
% \begin{align}
% 1 = \mathbf{g} \mathbf{W} = \left( \mathbf{g} \mathbf{e} \mathbf{e}^\top + \mathbf{g} (\mathbf{I} - \mathbf{e} \mathbf{e}^\top) \right) \mathbf{W} = (h_1 + h_2) \mathbf{W}
% \end{align}
% with $h_1 = \mathbf{g} \mathbf{e} \mathbf{e}^\top$ and $h_2 = \mathbf{g} (\mathbf{I} - \mathbf{e} \mathbf{e}^\top)$. 
% $h_1$ is a rank 1 matrix and $h_2 \mathbf{W} = 0$. We apply our standard truncation only to $h_2$ and after that set $\mathbf{g} = h_1 + T(h_2)$. 

\subsection{Stability of the proposed low-rank scheme}

It can be shown that the DLRA scheme proposed in \eqref{alg:DLRA} preserves the numerical stability of the full
conservative system presented in \eqref{eqs:KB-both}, which has been shown in Theorem \ref{theorem:stabilityMg}. The rewriting of equations \eqref{eqs:KB-both} into \eqref{eqs:KBrefDLRA-both}, the basis augmentations in \eqref{alg:KbasisAug} and \eqref{alg:LbasisAug}
and the implementation of the described truncation strategy are crucial for the proof. We
begin with the following definition.

\begin{definition}[Low-rank approximation of the fully discrete solution]
The DLRA approximation of the fully discrete solution $f$ of the linear Boltzmann-BGK equation at time $t_n$ is given by $\mathbf{f}^n = (f_{jk}^n) \in \mathbb{R}^{N_x \times N_v}$ with entries 
\begin{align*}
f_{jk}^n = \frac{1}{\sqrt{2\pi}} \rho_j^n \sum_{m,\ell=1}^{r} X_{jm}^n S_{m\ell}^n V_{k\ell}^n e^{-v_k^2/2}. 
\end{align*} 
\end{definition}

Note that in this notation we do not distinguish between the full solution $\mathbf{f}^n$ and its low-rank approximation $\mathbf{f}^n_r$ at time $t_n$. Then we can show that the DLRA scheme \eqref{alg:DLRA} is numerically stable.

We can then show that algorithm \eqref{alg:DLRA} is  numerically stable.

\begin{theorem}
Under the time step restriction $\max_k(|v_k|) \Delta t \leq \Delta x$, the fully discrete DLRA scheme \eqref{alg:DLRA} is numerically stable in the $\mathscr{H}$-norm, i.e
\begin{align*}
\Vert \mathbf{f}^{n+1}  \Vert_{\mathscr{H}}^2 \leq \Vert \mathbf{f}^n \Vert_{\mathscr{H}}^2.
\end{align*}
\end{theorem}
\begin{proof}
We begin with the $S$-step given in \eqref{alg:Sstep}, multiply it with $\doublehat X_{\alpha q}^{n+1}$ and $\doublehat V_{\beta p}^{n+1}$ and sum over $q$ and $p$. For simplicity of notation, we introduce the projections $P_{\alpha j}^{X}= \sum_{q=1}^{4r} \doublehat X_{\alpha q}^{n+1} \doublehat X_{jq}^{n+1}$ and $P_{k\beta}^{V} = \sum_{p=1}^{4r} \doublehat V_{kp}^{n+1} \doublehat V_{\beta p}^{n+1}$. This yields
\begin{align*}
g_{\alpha \beta}^{n+1} = \ \frac{1}{1+\sigma \Delta t} \Bigg(& \sum_{j=1}^{N_x} \sum_{k=1}^{N_v} \frac{\rho_j^n}{\rho_j^{n+1}} g_{jk}^n P_{\alpha j}^{X} P_{k\beta}^{V} - \Delta t \sum_{i,j=1}^{N_x} \sum_{k=1}^{N_v} \frac{1}{\rho_j^{n+1}} D_{ji}^x \rho_i^n g_{ik}^n v_k P_{\alpha j}^{X} P_{k\beta}^{V}\nonumber\\
&+ \Delta t \frac{\Delta x}{2} \sum_{i,j=1}^{N_x} \sum_{k=1}^{N_v} \frac{1}{\rho_j^{n+1}} D_{ji}^{xx} \rho_i^n g_{ik}^n |v_k| P_{\alpha j}^{X} P_{k\beta}^{V} + \sigma \Delta t \sum_{j=1}^{N_x} \sum_{k=1}^{N_v} P_{\alpha j}^{X} P_{k\beta}^{V} \Bigg).
\end{align*}
Multiplication with $\frac{2}{\sqrt{2\pi}} \left(\rho_\alpha^{n+1}\right)^2 g_{\alpha \beta}^{n+1} \omega_\beta e^{v_\beta^2/2} (1+\sigma \Delta t)$ and summation over $\alpha$ and $\beta$ leads to
\begin{align*}
&\ \frac{2}{\sqrt{2\pi}} \sum_{\alpha=1}^{N_x} \sum_{\beta=1}^{N_v} \left(\rho_\alpha^{n+1} g_{\alpha \beta}^{n+1} \right)^2 \omega_\beta e^{v_\beta^2/2} \left(1+\sigma \Delta t \right)\\
=& \ \frac{2}{\sqrt{2\pi}} \sum_{j,\alpha=1}^{N_x} P_{\alpha j}^{X} \frac{\rho_j^n}{\rho_j^{n+1}} \left(\rho_\alpha^{n+1}\right)^2 \sum_{k=1}^{N_v} g_{jk}^n \sum_{\beta=1}^{N_v} P_{k\beta}^{V} g_{\alpha \beta}^{n+1}  \omega_\beta e^{v_\beta^2/2}\\
&- \frac{2 \Delta t}{\sqrt{2\pi}} \sum_{j,\alpha=1}^{N_x} P_{\alpha j}^{X} \frac{1}{\rho_j^{n+1}} \left(\rho_\alpha^{n+1}\right)^2 \sum_{i=1}^{N_x} D_{ji}^x \rho_i^n \sum_{k=1}^{N_v} g_{ik}^n v_k \sum_{\beta=1}^{N_v} P_{k\beta}^{V} g_{\alpha \beta}^{n+1} \omega_\beta e^{v_\beta^2/2}\\
&+ \Delta t \frac{\Delta x}{\sqrt{2\pi}} \sum_{j,\alpha=1}^{N_x} P_{\alpha j}^{X} \frac{1}{\rho_j^{n+1}} \left(\rho_\alpha^{n+1}\right)^2 \sum_{i=1}^{N_x} D_{ji}^{xx} \rho_i^n \sum_{k=1}^{N_v} g_{ik}^n |v_k| \sum_{\beta=1}^{N_v} P_{k\beta}^{V} g_{\alpha \beta}^{n+1} \omega_\beta e^{v_\beta^2/2}\\
&+ \frac{2 \sigma \Delta t}{\sqrt{2\pi}} \sum_{j,\alpha=1}^{N_x} P_{\alpha j}^{X} \left(\rho_\alpha^{n+1} \right)^2 \sum_{\beta=1}^{N_v} g_{\alpha \beta}^{n+1} \omega_\beta e^{v_\beta^2/2}\sum_{k=1}^{N_v} P_{k\beta}^{V} .
\end{align*}
Using the basis augmentations given in \eqref{alg:KbasisAug} and \eqref{alg:LbasisAug}, we can deduce that the equalities
\begin{align*}
\sum_{\alpha=1}^{N_x} P_{\alpha j}^{X} g_{\alpha \beta}^{n+1} (\rho_\alpha^{n+1})^2 = g_{j \beta}^{n+1} (\rho_j^{n+1})^2 \quad \text{ and } \quad \sum_{\beta=1}^{N_v} P_{k\beta}^{V} g_{j\beta}^{n+1} \omega_\beta e^{v_\beta^2/2} = g_{jk}^{n+1} \omega_k e^{v_k^2/2}
\end{align*}
hold. We insert these relations and, to be consistent in notation, change the summation indices on the left-hand side from $\alpha$ to $j$ and from $\beta$ to $k$. This leads to
\begin{align*}
\frac{2}{\sqrt{2\pi}} \sum_{j=1}^{N_x} \sum_{k=1}^{N_v} \left(\rho_j^{n+1} g_{jk}^{n+1} \right)^2 \omega_k e^{v_k^2/2} \left(1+\sigma \Delta t \right)=& \ \frac{2}{\sqrt{2\pi}} \sum_{j=1}^{N_x} \sum_{k=1}^{N_v} \rho_j^n g_{jk}^n \rho_j^{n+1} g_{jk}^{n+1} \omega_k e^{v_k^2/2}\\
&- \frac{2 \Delta t}{\sqrt{2\pi}} \sum_{i,j=1}^{N_x} \sum_{k=1}^{N_v} \rho_j^{n+1} g_{jk}^{n+1} D_{ji}^x \rho_i^n g_{ik}^n v_k \omega_k e^{v_k^2/2}\\
&+ \Delta t \frac{\Delta x}{\sqrt{2\pi}} \sum_{i,j=1}^{N_x} \sum_{k=1}^{N_v} \rho_j^{n+1} g_{jk}^{n+1} D_{ji}^{xx} \rho_i^n g_{ik}^n |v_k| \omega_k e^{v_k^2/2}\\
&+ \frac{2 \sigma \Delta t}{\sqrt{2\pi}} \sum_{j=1}^{N_x} \sum_{k=1}^{N_v}  \left(\rho_j^{n+1}\right)^2 g_{jk}^{n+1} \omega_k e^{v_k^2/2}.
\end{align*}
% We bring the term $\frac{2 \sigma \Delta t}{\sqrt{2\pi}} \sum_{j=1}^{N_x} \sum_{k=1}^{N_v} \left(\rho_j^{n+1} g_{jk}^{n+1} \right)^2 \omega_k e^{v_k^2/2}$ to the right-hand side and use the relation 
% \begin{align*}
% 2 \rho_j^n g_{jk}^n \rho_j^{n+1} g_{jk}^{n+1}  = \left( \rho_j^{n+1} g_{jk}^{n+1} \right)^2 +  \left( \rho_j^n g_{jk}^n \right)^2 - \left( \rho_j^{n+1} g_{jk}^{n+1} - \rho_j^n g_{jk}^n \right)^2.
% \end{align*}
% With the notations $f_{jk}^{n,n+1} = \frac{1}{\sqrt{2\pi}} e^{-v_k^2/2} \rho_j^{n,n+1} g_{jk}^{n,n+1}$ and $M_{jk}^{n+1} = \frac{1}{\sqrt{2\pi}}  \rho_j^{n+1} e^{-v_k^2/2}$, we then obtain
% \begin{align*}
% \Vert \mathbf{f}^{n+1} \Vert_{\mathscr{H}}^2 =& \ \Vert \mathbf{f}^n \Vert_{\mathscr{H}}^2 - \sqrt{2\pi} \sum_{j=1}^{N_x} \sum_{k=1}^{N_v}(f_{jk}^{n+1} - f_{jk}^n)^2 \omega_k e^{3v_k^2/2}\\
% &-2 \sqrt{2\pi} \Delta t \sum_{i,j=1}^{N_x} \sum_{k=1}^{N_v} f_{jk}^{n+1} D_{ji}^x f_{ik}^n v_k \omega_k e^{3v_k^2/2} + 2 \sqrt{2\pi} \Delta t \sum_{i,j=1}^{N_x} \sum_{k=1}^{N_v} f_{jk}^{n+1} D_{ji}^{xx} f_{ik}^n |v_k| \omega_k e^{3v_k^2/2}\\
% &+2 \sqrt{2\pi} \sigma \Delta t \sum_{j=1}^{N_x} \sum_{k=1}^{N_v} f_{jk}^{n+1} \left(M_{jk}^{n+1} -f_{jk}^{n+1}\right) \omega_k e^{3v_k^2/2},
% \end{align*}
Analogously to the proof of Theorem \ref{theorem:stabilityMg} and, as the truncation step is designed not to alter these expressions, we can conclude that the proposed DLRA scheme decreases the $\mathscr{H}$-norm and hence is numerically stable under the time step restriction $\max_k(|v_k|) \Delta t \leq \Delta x$. 
\end{proof}

\begin{remark}
In case of numerical instabilities occurring in practical applications due to very small $\rho >0$, standard numerical methods such as slope limiters can be added.
\end{remark}

\section{Numerical results}\label{sec6:NumericalResults}

To validate the theoretical considerations and the proposed DLRA scheme, numerical results for different test examples in 1D as well as in 2D are given in this section. 

\subsection{1D Plane source}\label{sec:1DPlanesource}

We begin with a one-dimensional analogue to the plane source problem, which is a common test case for radiative transfer \cite{baumann2024, ceruti2022rank, kusch2023stability,peng2020-2D}, and compare the solution of the full equations \eqref{eqs:KB-both} to the solution obtained by the DLRA scheme given in \eqref{alg:DLRA}. The spatial domain shall be set to $D=[-10,10]$. As initial conditions we choose the density $\rho$ to be a cutoff Gaussian
\begin{align*}
\rho(t=0,x) = \max \left( 10^{-4}, \frac{1}{\sqrt{2\pi \sigma_{\text{IC}}^2}} \exp{\left( - \frac{x^2}{2 \sigma_{\text{IC}}^2} \right) } \right)
\end{align*}
with constant deviation $\sigma_{\text{IC}} = 0.3$. The initial distribution function $g$ is assumed to be constant in space and velocity and we prescribe $g(t=0,x,v)=1$. We consider relatively large collisionality by choosing $\sigma=10$. Computations are started with an initial rank of $r=20$. As computational parameters we use $N_x = 1000$ grid points in the spatial as well as $N_v = 500$ grid points in the velocity domain. Due to this choice, we obtain $\max_k(|v_k|) \approx 31.05$, which is adjusted to the next larger integer value. The time step size is determined by $\Delta t = \text{CFL} \cdot \frac{\Delta x}{32}$ with $\text{CFL}=0.99$, according to the corresponding CFL condition. 

Practical implementations show that the basis augmentations to rank $4r$ in in \eqref{alg:KbasisAug} and \eqref{alg:LbasisAug}, which are needed for the theoretical proof of the numerical stability, may not be necessary for numerical examples and that the standard basis augmentations to rank $2r$ provide similar solutions while being significantly faster. For this reason, we propose to leave out the basis augmentations \eqref{alg:KbasisAug} and \eqref{alg:LbasisAug} in practical applications. In this case, all quantities with double hats related to rank $4r$ reduce to quantities of rank $2r$ with one single hat. The simplified scheme with rank $2r$ is also visualized (in brackets) in the flowchart of Figure \ref{fig:flowchart}. 

\begin{figure}[h!]
    \centering
    \includegraphics[width = 0.95\linewidth]{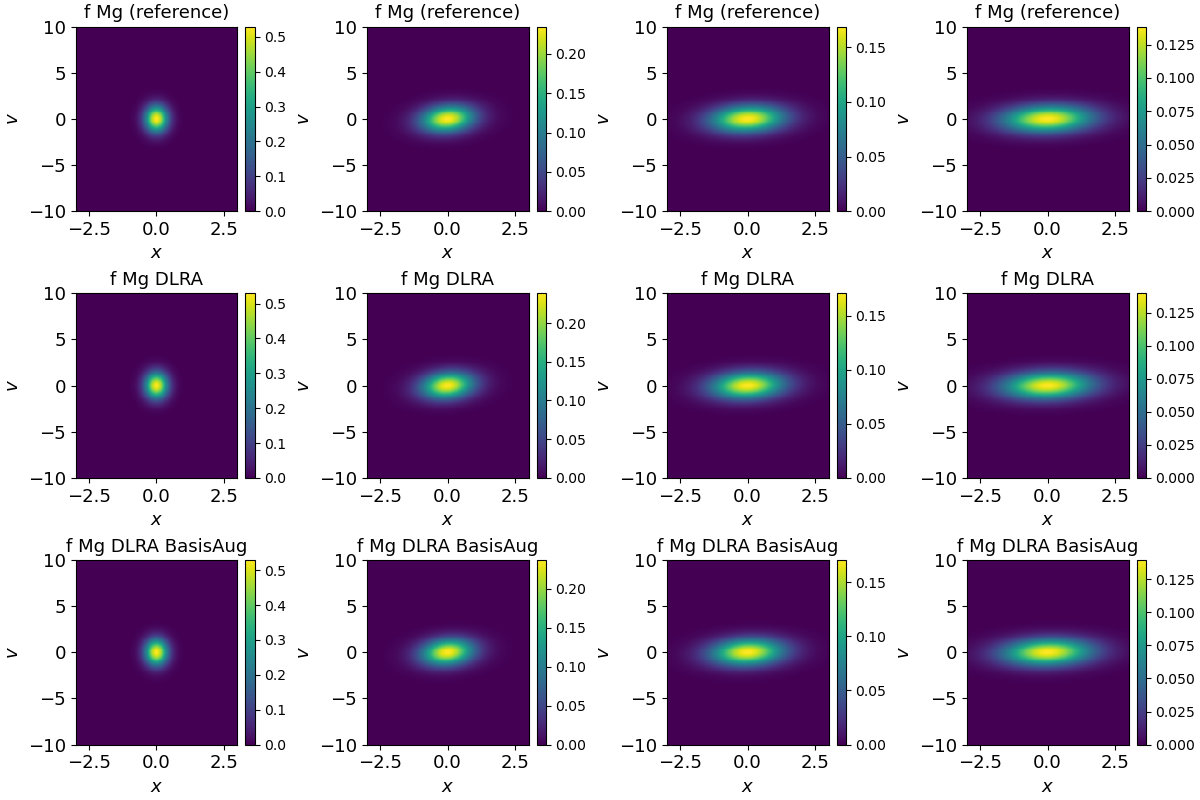}
   \caption{Numerical results for the solution $f(t,x,v)$ of the 1D plane source analogue at time $t=0$ (first column), $t=2$ (second column), $t=4$ (third column), and $t=6$ (fourth column), computed with the full solver (Mg (reference)) (first row), the reduced DLRA scheme (Mg DLRA) (second row) and the basis augmented DLRA scheme (Mg DLRA BasisAug) (third row).}
    \label{fig:1DPlanesource-f}
\end{figure}

\begin{figure}[h!]
    \centering
    \includegraphics[width = 0.9\linewidth]{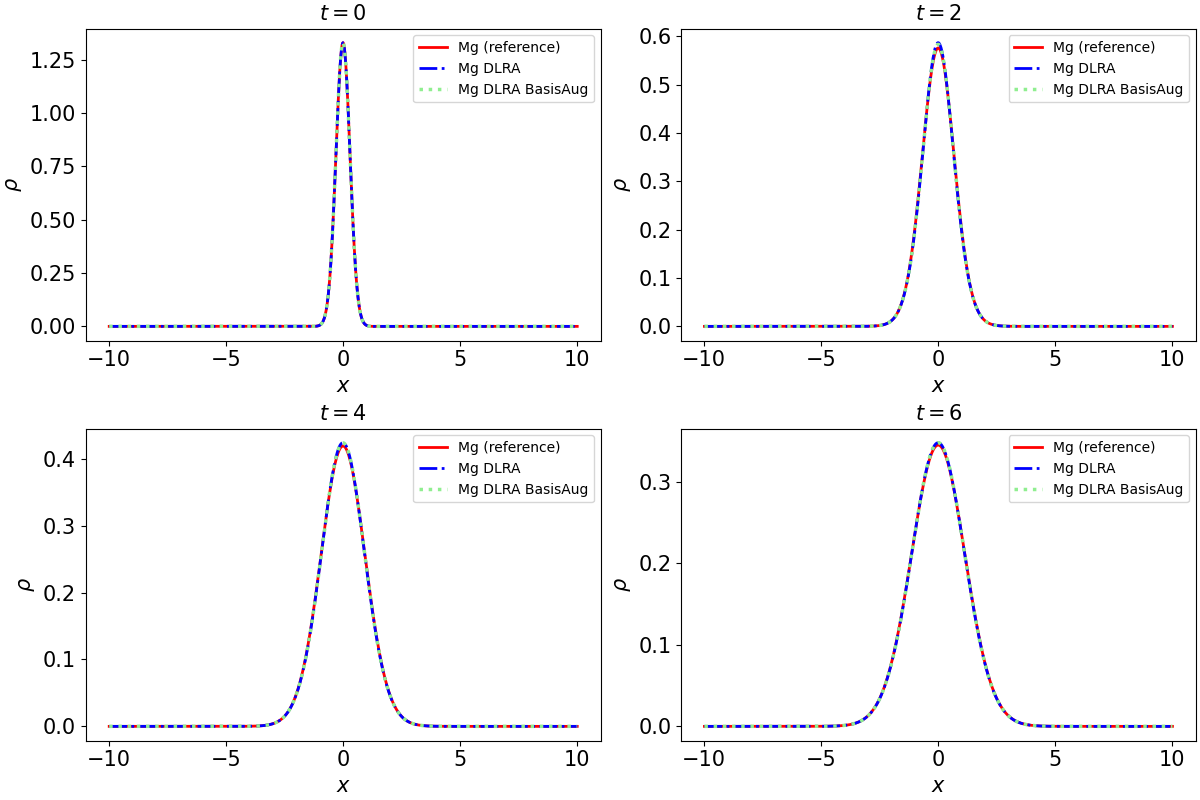}
   \caption{Numerical results for the density $\rho(t, x)$ of the 1D plane source analogue at time $t=0$, $t=2$, $t=4$, and $t=6$, computed with the full solver (Mg (reference)), the reduced DLRA scheme (Mg DLRA) and the basis augmented DLRA scheme (Mg DLRA BasisAug).}
    \label{fig:1DPlanesource-rho}
\end{figure}

\begin{figure}[h!]
    \centering
    \includegraphics[width = 1.0\linewidth]{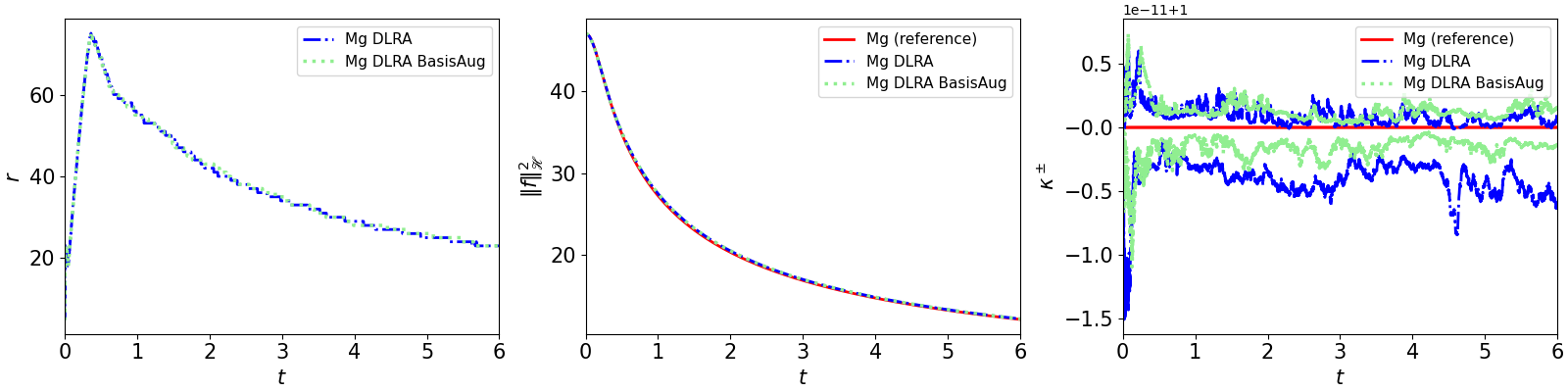}
   \caption{Left: Evolution of the rank in time for the 1D plane source analogue for the reduced DLRA scheme (Mg DLRA) and the basis augmented DLRA scheme (Mg DLRA BasisAug). Middle: Evolution of the $\mathscr{H}$-norm in time for the full solver (Mg (reference)), the reduced DLRA scheme (Mg DLRA) and the basis augmented DLRA scheme (Mg DLRA BasisAug). Right: Evolution of $\kappa^\pm$ in time for the full solver (Mg (reference)), the reduced DLRA scheme (Mg DLRA) and the basis augmented DLRA scheme (Mg DLRA BasisAug). The line corresponding to the full system has the constant value $1$.}
    \label{fig:1DPlanesource-rank}
\end{figure}

In Figure \ref{fig:1DPlanesource-f}, we compare the results for the full solution $f(t,x,v)$, computed with the full solver (Mg (reference)), the DLRA scheme with rank $2r$ (Mg DLRA) and the basis augmented DLRA scheme with rank $4r$ (Mg DLRA BasisAug) at different times up to $t_{\text{End}}=6$. We observe that the reduced as well as the augmented DLRA algorithm capture the main characteristics of the full reference Mg system.
This is also true for the computational results for the density $\rho(t,x)$, displayed in Figure \ref{fig:1DPlanesource-rho}.
Figure \ref{fig:1DPlanesource-rank} shows the evolution of the rank, which for a chosen tolerance parameter of $\vartheta = 10^{-5} \Vert \mathbf{\Sigma} \Vert_2$ increases up to $r=75$ before it significantly reduces over time. It confirms that even though the initial condition is chosen to be of low-rank, the evolution of the solution in time can leave the low-rank regime. Note that the evolution of the rank for the reduced as well as for the basis augmented algorithm show good agreement as the new rank is displayed after the corresponding truncation step. Further, the behavior of the norm $\Vert \mathbf{f} \Vert_\mathscr{H}^2$ is depicted. As expected, it decreases smoothly over time for all considered systems. Additionally, we display the quantities $\kappa^+ := \max_j \left(\frac{1}{\sqrt{2\pi}} \sum_{k=1}^{N_v} g_{jk} \omega_k e^{v_k^2/2}\right)$ and $\kappa^- := \min_j \left(\frac{1}{\sqrt{2\pi}} \sum_{k=1}^{N_v} g_{jk} \omega_k e^{v_k^2/2}\right)$. According to Lemma \ref{lemma:fullydiscrete-rho} it is essential that they are equal to $1$, which for the low-rank schemes is ensured by the adjusted truncation step. It can be observed that this property is fulfilled up to order $\mathcal{O}\left(10^{-11}\right)$.

\subsection{2D Plane source}\label{sec6.2:2D Plane Source}

In higher dimensions, the computational advantages of DLRA schemes are enhanced. For this reason, we give some two-dimensional test examples starting with the two-dimensional version of the plane source problem considered in the previous section. The corresponding two-dimensional set of equations becomes
\begin{align*}
\partial_t g(t,\mathbf{x},\mathbf{v}) &= - \frac{\mathbf{v}}{\rho(t,\mathbf{x})} \cdot \nabla_\mathbf{x} \left(\rho(t,\mathbf{x}) g(t,\mathbf{x},\mathbf{v})\right) + \sigma \left(1-g(t,\mathbf{x},\mathbf{v})\right) - \frac{g(t,\mathbf{x},\mathbf{v})}{\rho(t,\mathbf{x})} \partial_t \rho(t,\mathbf{x}),\\
\partial_t \rho(t,\mathbf{x}) &= - \frac{1}{2\pi} \nabla_\mathbf{x} \cdot \int \rho(t,\mathbf{x}) g(t,\mathbf{x},\mathbf{v}) \mathbf{v} e^{-|\mathbf{v}|^2/2} \mathrm{d}\mathbf{v},
\end{align*}
where $\mathbf{x}=(x_1,x_2) \in D \subset \mathbb{R}^2$ and $\mathbf{v}=(v_1,v_2) \in \mathbb{R}^2$. We compare the solution of this system with the solution of the DLRA scheme for which the extension to the two-dimensional setting is straightforward. For this test example we choose the spatial domain $D=[-3,3] \times [-3,3]$ and prescribe the initial condition for the density by 
\begin{align*}
\rho(t=0,\mathbf{x}) = \frac{1}{4\pi} \max\left( 10^{-1}, \frac{10^2}{4\pi \sigma_{\text{IC}}^2} \exp{\left( - \frac{|\mathbf{x}|^2}{4 \sigma_{\text{IC}}^2} \right)} \right) 
\end{align*}
with constant deviation $\sigma_{\text{IC}}=0.3$. The function $g$ is set to $g(t=0,\mathbf{x},\mathbf{v})=1$ and the collision coefficient to $\sigma=100$. We allow an initial rank of $r=20$. Computations are performed on a spatial grid with $N_{x_1}=128$ grid points in $x_1$ and $N_{x_2}=128$ grid points in $x_2$. For the velocity grid we choose $N_{v_1}=32$ grid points in $v_1$ and $N_{v_2}=32$ grid points in $v_2$. Due to this choice, we obtain $\max_k \left(|\mathbf{v}_k| \right) \approx 10.08,$ which is adjusted to the next larger integer value. The time step size is determined by $\Delta t = \text{CFL} \cdot \frac{\Delta x}{11}$ with a CFL number of $\text{CFL}=0.7$. 

\begin{figure}[t]
    \centering
    \includegraphics[width = 1.0\linewidth]{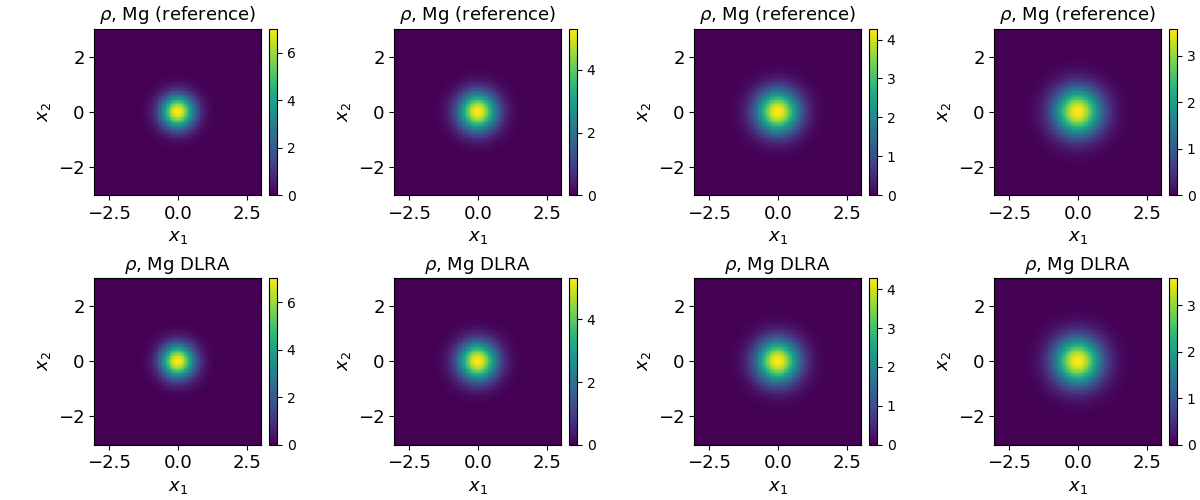}

   \caption{Numerical results for the density $\rho(t, \mathbf{x})$ of the 2D plane source analogue at time $t=0$ (first column), $t=1$ (second column), $t=2$ (third column, and $t=3$ (fourth column), computed with the full solver (Mg (reference)) (first row) and the reduced DLRA scheme (Mg DLRA) (second row).}
    \label{fig:2DPlanesource-rho}
\end{figure}

\begin{figure}[t]
    \centering
    \includegraphics[width = 1.0\linewidth]{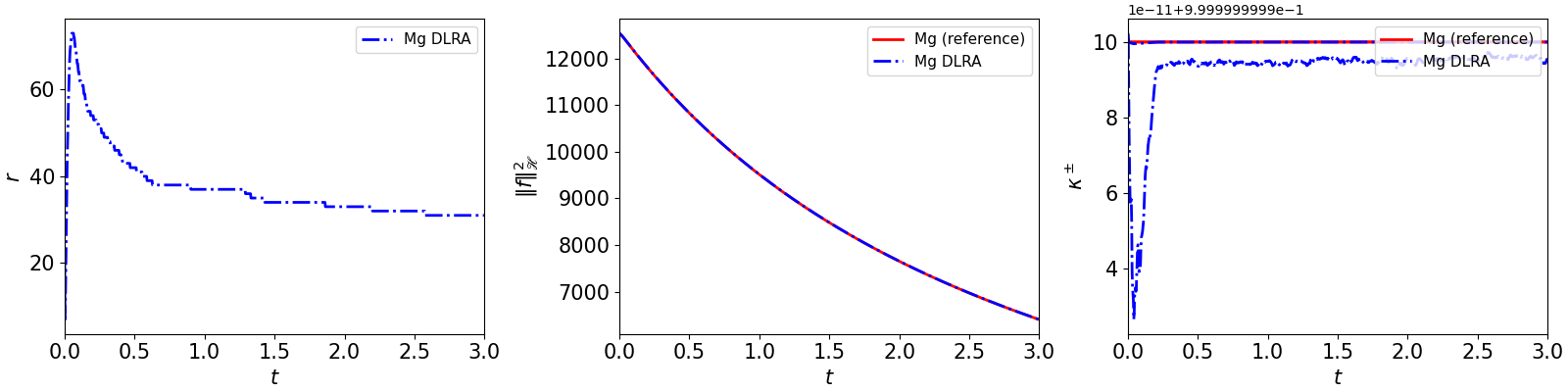}
   \caption{Left: Evolution of the rank in time for the 2D plane source analogue for the reduced DLRA scheme (Mg DLRA). Middle: Evolution of the $\mathscr{H}$-norm in time for the full solver (Mg (reference)) and the reduced DLRA scheme (Mg DLRA). Right: Evolution of $\kappa^\pm$ in time for the full solver (Mg (reference)) and the reduced DLRA scheme (Mg DLRA).}
    \label{fig:2DPlanesource-rank}
\end{figure}

Figure \ref{fig:2DPlanesource-rho} compares the density $\rho(t,x)$ at different times up to $t_{\text{End}}=3.0$, computed with the full solver (Mg (reference)) and the reduced DLRA scheme with rank $2r$ (Mg DLRA). Note that we refrain from computations with the basis augmented $4r$ scheme as in two space and velocity dimensions this would lead to extremely increased computational costs. We observe that the solution of the DLRA scheme matches the solution of the full system.
To determine the evolution of the rank, we use a tolerance parameter of $\vartheta = 10^{-5} \Vert \mathbf{\Sigma} \Vert_2$. In Figure \ref{fig:2DPlanesource-rank}, we observe an increasing up to $r=73$ before it decreases continuously over time. Further, the evolution of the norm $\Vert \mathbf{f} \Vert_\mathscr{H}^2$ in time is displayed. It decreases smoothly over time for all considered systems. In addition, we plot the quantities $\kappa^+ := \max_j \left(\frac{1}{2\pi} \sum_{k=1}^{N_{v_1}} \sum_{\ell=1}^{N_{v_2}} g(t,\mathbf{x}_j,v_k^1, v_\ell^2) \omega_k \omega_\ell e^{\left(v_k^1 \right)^2/2} e^{\left(v_\ell^2 \right)^2/2}\right)$ and $\kappa^- := \min_j \left(\frac{1}{2\pi} \sum_{k=1}^{N_{v_1}} \sum_{\ell=1}^{N_{v_2}} g(t,\mathbf{x}_j,v_k^1, v_\ell^2) \omega_k \omega_\ell e^{\left(v_k^1 \right)^2/2} e^{\left(v_\ell^2 \right)^2/2}\right)$, which are equal to $1$ up to order $\mathcal{O}\left(10^{-10}\right)$.
For this setup, the running time of the DLRA scheme compared to the full solver is clearly faster. It reduces by a factor of approximately 10 from 2315 seconds to 235 seconds, confirming the computational advantages of the DLRA scheme.

\subsection{2D Beam}

% User defined $\mathbf{v}_{\mathrm{beam}} = \frac{1}{\sqrt{2}}\begin{pmatrix}
%     -1 \\ -1
% \end{pmatrix}$ as an example.
% \begin{align*}
%     \tilde g_{jk} = \exp\left(-\frac{1}{\sigma^2}\Vert \mathbf{v}_k - \mathbf{v}_{\mathrm{beam}} \Vert^2 \right) \\
%     g_{jk} = \tilde g_{jk} \cdot \frac{1}{\frac{1}{\sqrt{2\pi}}\sum_{\ell} \omega_\ell e^{v_\ell^2/2} \tilde g_{j\ell}}
% \end{align*}
% Then,
% \begin{align*}
%     \frac{1}{\sqrt{2\pi}}\sum_k \omega_k e^{v_k^2/2} g_{jk} = \frac{1}{\sqrt{2\pi}} \sum_k \omega_k e^{v_k^2/2} \tilde g_{jk} \cdot \frac{1}{\frac{1}{\sqrt{2\pi}} \sum_{\ell} \tilde \omega_\ell e^{v_\ell^2/2} \tilde g_{j\ell}} = 1
% \end{align*}

As a second two-dimensional test example we consider a beam in the spatial domain $D=[-5,5] \times [-5,5]$ starting at $(0,0)$ in the middle of the spatial plane and moving to the bottom left. The initial values are given by
\begin{align*}
\rho(t=0,\mathbf{x}) &= \frac{1}{4\pi} \max \left( 10^{-1}, \frac{10^2}{4\pi \sigma_{\text{IC},\rho}^2} \exp{\left( - \frac{|\mathbf{x}|^2}{4 \sigma_{\text{IC},\rho}^2} \right) } \right),\\
g(t=0, \mathbf{x}, \mathbf{v}) &= \frac{C}{4\pi} \max\left( 10^{-14}, \frac{10^6}{4\pi \sigma_{\text{IC},g}^2} \exp\left(-\frac{\left| \mathbf{v} - \mathbf{v}_{\mathrm{beam}} \right|^2}{4 \sigma_{\text{IC},g}^2} \right)\right),
\end{align*}
where $\sigma_{\text{IC},\rho}=0.2$, $\sigma_{\text{IC},g}=0.01$ and $C$ is a normalization constant such that $\frac{1}{2\pi} \int g(t=0, \mathbf{x}, \mathbf{v}) e^{-|\mathbf{v}|^2/2} \mathrm{d}\mathbf{v}= \mathbf{1}$ holds. The beam velocity $\mathbf{v}_{\mathrm{beam}}$ is set to $\mathbf{v}_{\mathrm{beam}} = \begin{pmatrix}
    -1 \\ -1
\end{pmatrix}$ and the collisionality to $\sigma = 1.5$. All other initial settings and computational parameters remain unchanged from the previous test example given in Section \ref{sec6.2:2D Plane Source}.

\begin{figure}[t]
    \centering
    \includegraphics[width = 1.0\linewidth]{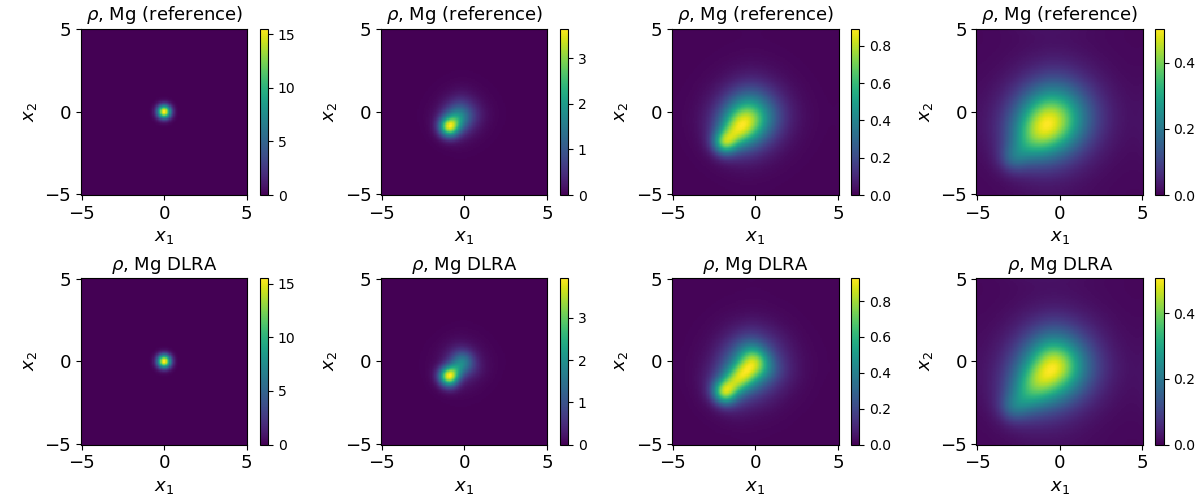}
   \caption{Numerical results for the density $\rho(t, \mathbf{x})$ of the 2D beam test problem at time $t=0$ (first column), $t=1$ (second column), $t=2$ (third column, and $t=3$ (fourth column), computed with the full solver (Mg reference)) (first row) and the reduced DLRA scheme (Mg DLRA) (second row).}
    \label{fig:2DBeam-rho}
\end{figure}

\begin{figure}[t]
    \centering
    \includegraphics[width = 1.0\linewidth]{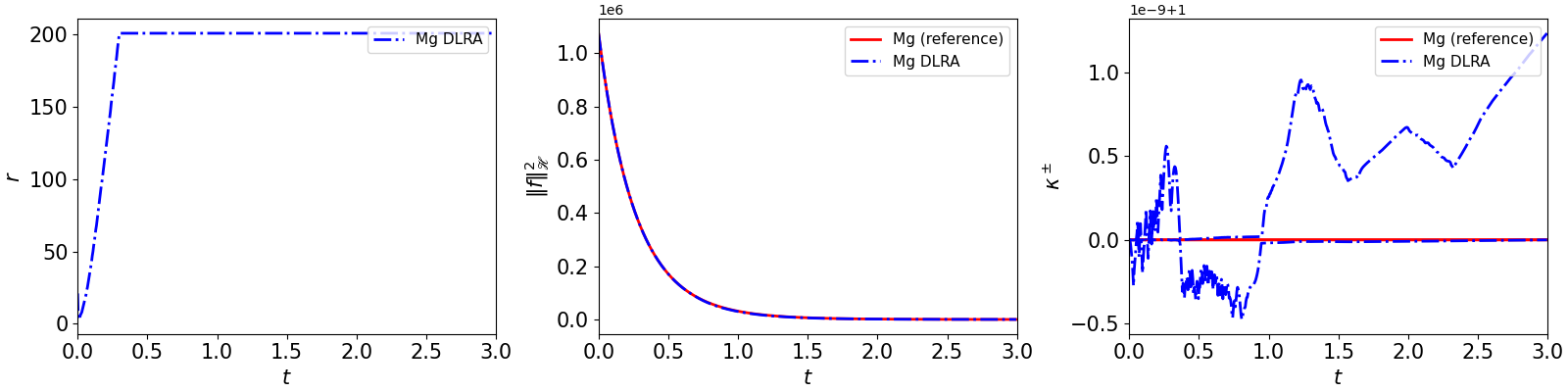}
   \caption{Left: Evolution of the rank in time for the 2D beam test problem for the reduced DLRA scheme (Mg DLRA). The rank increases up to the maximal allowed value of $r=200$. Middle: Evolution of the $\mathscr{H}$-norm in time for the full solver (Mg (reference)) and the reduced DLRA scheme (Mg DLRA). Right: Evolution of $\kappa^\pm$ in time for the full solver (Mg (reference)) and the reduced DLRA scheme (Mg DLRA). The red line has the constant value $1$. The deviation of the DLRA scheme from $1$ is of order $\mathcal{O}\left(10^{-12}\right)$.}
    \label{fig:2DBeam-rank}
\end{figure}

Figure \ref{fig:2DBeam-rho} compares the density $\rho(t,x)$ at different times up to $t_{\text{End}}=3.0$, computed with the full solver (Mg (reference)) and the reduced DLRA scheme with rank $2r$ (Mg DLRA). At all displayed time steps the DLRA solution resembles the solution of the full system.
In Figure \ref{fig:2DBeam-rank}, the evolution of the rank in time is shown. We use a tolerance parameter of $\vartheta = 10^{-4} \Vert \mathbf{\Sigma} \Vert_2$ and allow a maximal rank of $200$. Due to the choice of $\sigma$ the solution of the problem is not low-rank. For this reason, we observe an increasing of the rank from rank $1$ for the initially chosen configuration up to the maximal allowed value of $r=200$. Also, the evolution of the norm $\Vert \mathbf{f} \Vert_\mathscr{H}^2$ in time is depicted. It decreases continuously, matching our theoretical considerations. Further, it can be observed that the quantities $\kappa^+$ and $\kappa^-$ defined in Section \ref{sec6.2:2D Plane Source} are equal to $1$ up to order $\mathcal{O}\left(10^{-9}\right)$.
Due to the high rank, the computational benefits of the DLRA scheme are diminished compared to the previous test case. It reduces by a factor of approximately $1.5$ from $1220$ seconds to $845$ seconds. This example illustrates the relation between the choice of $\sigma$ and the low-rank structure of the solution. It is expected that for larger values of $\sigma$ the solution becomes low-rank and hence the computational benefits of the DLRA scheme are enhanced.

\section{Conclusion and outlook}\label{sec7:Outlook}

We have derived a multiplicative DLRA discretization for the linear Boltzmann-BGK problem that in contrast to another presented naive discretization is numerically stable. To show this, we have conducted a stability analysis leading to a concrete hyperbolic CFL condition. In addition, numerical examples in 1D and 2D confirm the stability, accuracy and efficiency of the proposed DLRA scheme.\\
The insights gained from this article can be helpful for future work as the employed multiplicative splitting is attached to the investigation of more complicated equations, e.g. the non-linear Boltzmann-BGK equation treated in \cite{einkemmerhuying2021}. However, a direct transition of knowledge concerning the theoretical stability analysis is hardly possible as for the non-linear case most of the theoretical concepts applied here (such as von Neumann stability analysis) are not available. Concerning the extension of the design of the numerical scheme from the linear problem to the non-linear setting, we expect a relatively straightforward extension to the weakly compressible case to be possible. In this case we can expand the Maxwellian distribution  function $M = \frac{\rho}{\sqrt{2\pi}} \exp{\left(- \frac{\left|v -u\right|^2}{2 \theta}\right)}$ in terms of the small (compared to the speed of sound) macroscopic velocity $u$. Also the extension to the general non-linear case is possible but a discretization of the conservative form of the equations, i.e. by not splitting up the term $\partial_x\left(Mg\right)$, cannot be efficiently implemented as the Maxwellian $M$ is generally not low-rank.

\section*{Acknowledgements}

The authors thank Marlies Pirner for helpful discussions on the linear Boltzmann-BGK equation. Lena Baumann acknowledges support by the Würzburg Mathematics Center for Communication and Interaction (WMCCI) as well as the Stiftung der Deutschen Wirtschaft.

\newpage
\bibliographystyle{abbrv}
\bibliography{references}

\end{document}